
\def\LaTeX{%
  \let\Begin\begin
  \let\End\end
  \let\salta\relax
  \let\finqui\relax
  \let\futuro\relax}

\def\UK{\def\our{our}\let\sz s}
\def\USA{\def\our{or}\let\sz z}

\UK



\LaTeX

\USA


\salta

\documentclass[twoside,12pt]{article}
\setlength{\textheight}{24cm}
\setlength{\textwidth}{16cm}
\setlength{\oddsidemargin}{2mm}
\setlength{\evensidemargin}{2mm}
\setlength{\topmargin}{-15mm}
\parskip2mm


\usepackage[usenames,dvipsnames]{color}
\usepackage{amsmath}
\usepackage{amsthm}
\usepackage{amssymb}
\usepackage[mathcal]{euscript}
%
%
\usepackage{cite}
%
%
%


\definecolor{viola}{rgb}{0.3,0,0.7}
\definecolor{ciclamino}{rgb}{0.5,0,0.5}
\definecolor{rosso}{rgb}{0.8,0,0}

\def\lastrev #1{{\color{rosso}#1}}

\def\gianni #1{#1}
\def\pier #1{#1}
\def\juerg #1{#1}
\def\lastrev #1{{#1}}



\bibliographystyle{plain}


%

\finqui

\def\Beq{\Begin{equation}}
\def\Eeq{\End{equation}}
\def\Bsist{\Begin{eqnarray}}
\def\Esist{\End{eqnarray}}

\def\Bthm{\Begin{theorem}}
\def\Ethm{\End{theorem}}
\def\Blem{\Begin{lemma}}
\def\Elem{\End{lemma}}
\def\Bprop{\Begin{proposition}}
\def\Eprop{\End{proposition}}
\def\Bcor{\Begin{corollary}}
\def\Ecor{\End{corollary}}
\def\Brem{\Begin{remark}\rm}
\def\Erem{\End{remark}}

\def\Bdim{\Begin{proof}}
\def\Edim{\End{proof}}
\def\Bcenter{\Begin{center}}
\def\Ecenter{\End{center}}
\let\non\nonumber




\def\step #1 \par{\medskip\noindent{\bf #1.}\quad}


\def\Lip{Lip\-schitz}
\def\Holder{H\"older}
\def\Frechet{Fr\'echet}
\def\aand{\quad\hbox{and}\quad}

\def\lhs{left-hand side}
\def\rhs{right-hand side}
\def\sfw{straightforward}


\def\generaliz{generali\sz}


\def\multibold #1{\def\arg{#1}%
  \ifx\arg\pto \let\next\relax
  \else
  \def\next{\expandafter
    \def\csname #1#1#1\endcsname{{\bf #1}}%
    \multibold}%
  \fi \next}

\def\pto{.}

\def\multical #1{\def\arg{#1}%
  \ifx\arg\pto \let\next\relax
  \else
  \def\next{\expandafter
    \def\csname cal#1\endcsname{{\cal #1}}%
    \multical}%
  \fi \next}


\def\multimathop #1 {\def\arg{#1}%
  \ifx\arg\pto \let\next\relax
  \else
  \def\next{\expandafter
    \def\csname #1\endcsname{\mathop{\rm #1}\nolimits}%
    \multimathop}%
  \fi \next}

\multibold
qwertyuiopasdfghjklzxcvbnmQWERTYUIOPASDFGHJKLZXCVBNM.

\multical
QWERTYUIOPASDFGHJKLZXCVBNM.

\multimathop
ad dist div dom meas sign supp .


\def\accorpa #1#2{\eqref{#1}--\eqref{#2}}
\def\Accorpa #1#2 #3 {\gdef #1{\eqref{#2}--\eqref{#3}}%
  \wlog{}\wlog{\string #1 -> #2 - #3}\wlog{}}


\def\infess{\mathop{\rm inf\,ess}}
\def\supess{\mathop{\rm sup\,ess}}

\def\neto{\mathrel{{\scriptscriptstyle\nearrow}}}
\def\seto{\mathrel{{\scriptscriptstyle\searrow}}}

\def\graffe #1{\mathopen\{#1\mathclose\}}

\def\<#1>{\mathopen\langle #1\mathclose\rangle}
\def\norma #1{\mathopen \| #1\mathclose \|}

\def\iot {\int_0^t}
\def\ioT {\int_0^T}
\def\intQt{\int_{Q_t}}
\def\intQ{\int_Q}
\def\iO{\int_\Omega}
\def\iG{\int_\Gamma}
\def\intS{\int_\Sigma}
\def\intSt{\int_{\Sigma_t}}

\def\dt{\partial_t}
\def\dn{\partial_n}

\def\cpto{\,\cdot\,}

\def\checkmmode #1{\relax\ifmmode\hbox{#1}\else{#1}\fi}
\def\aeO{\checkmmode{a.e.\ in~$\Omega$}}
\def\aeQ{\checkmmode{a.e.\ in~$Q$}}

\def\aeS{\checkmmode{a.e.\ on~$\Sigma$}}
\def\aet{\checkmmode{a.e.\ in~$(0,T)$}}

\def\aat{\checkmmode{for a.a.~$t\in(0,T)$}}


\def\erre{{\mathbb{R}}}




\def\genspazio #1#2#3#4#5{#1^{#2}(#5,#4;#3)}
\def\spazio #1#2#3{\genspazio {#1}{#2}{#3}T0}

\def\L {\spazio L}
\def\H {\spazio H}
\def\W {\spazio W}

\def\C #1#2{C^{#1}([0,T];#2)}


\def\Lx #1{L^{#1}(\Omega)}
\def\Hx #1{H^{#1}(\Omega)}

\def\LxG #1{L^{#1}(\Gamma)}
\def\HxG #1{H^{#1}(\Gamma)}

\def\LQ #1{L^{#1}(Q)}
\def\LS #1{L^{#1}(\Sigma)}

\def\Ldue{\Lx 2}

\def\Huno{\Hx 1}
\def\Hdue{\Hx 2}
\def\Hunoz{{H^1_0(\Omega)}}
\def\HunoG{\HxG 1}
\def\HdueG{\HxG 2}

\def\LdueG{\LxG 2}


\def\LQ #1{L^{#1}(Q)}


\let\theta\vartheta

\let\phi\varphi
\let\lam\lambda

\let\TeXchi\chi                         
\newbox\chibox
\setbox0 \hbox{\mathsurround0pt $\TeXchi$}
\setbox\chibox \hbox{\raise\dp0 \box 0 }
\def\chi{\copy\chibox}



\def\suG{{\vrule height 5pt depth 4pt\,}_\Gamma}

\def\fG{f_\Gamma}
\def\yG{y_\Gamma}
\def\uG{u_\Gamma}
\def\vG{v_\Gamma}
\def\hG{h_\Gamma}
\def\xiG{\xi_\Gamma}
\def\qG{q_\Gamma}
\def\uGmin{u_{\Gamma,{\rm min}}}
\def\uGmax{u_{\Gamma,{\rm max}}}

\def\yz{y_0}

\def\Mz{M_0}
\def\ustar{u_*}
\def\vstar{v_*}
\def\vz{v_0}

\def\zT{\zzz_T}

\def\bQ{b_Q}
\def\bS{b_\Sigma}
\def\bO{b_\Omega}
\def\bG{b_\Gamma}
\def\bz{b_0}
\def\Az{\calA_0}

\def\rmin{r_-}
\def\rmax{r_+}

\def\zQ{z_Q}
\def\zS{z_\Sigma}
\def\zO{z_\Omega}
\def\zG{z_\Gamma}
\def\phQ{\phi_Q}
\def\phS{\phi_\Sigma}
\def\phO{\phi_\Omega}
\def\phG{\phi_\Gamma}

\def\Uad{\calU_{\pier{\ad}}}
\def\uopt{\overline u_\Gamma}
\def\yopt{\overline y}
\def\yGopt{\overline y_\Gamma}
\def\wopt{\overline w}

\def\redJ{\widetilde\calJ}

\def\unG{u_{\Gamma,n}}
\def\yn{y_n}
\def\ynG{y_{\Gamma,n}}
\def\wn{w_n}

\def\yh{y^h}
\def\yhG{y^h_\Gamma}
\def\wh{w^h}
\def\qh{q^h}
\def\qhG{q^h_\Gamma}
\def\zh{z^h}
\def\psiG{\psi_\Gamma}

\def\pO{p^\Omega}

\def\uO{u^\Omega}
\def\vO{v^\Omega}
\def\phOO{(\phO)^\Omega}

\def\VO{\calV_\Omega}
\def\HO{\calH_\Omega}
\def\RO{\calR_\Omega}
\def\VOp{\VO^*}

\def\VG{V_\Gamma}
\def\HG{H_\Gamma}
\def\Vp{V^*}

\def\normaV #1{\norma{#1}_V}
\def\normaH #1{\norma{#1}_H}

\def\normaHG #1{\norma{#1}_{\HG}}
\def\normaVp #1{\norma{#1}_*}

\def\nablaG{\nabla_{\!\Gamma}}
\def\DeltaG{\Delta_\Gamma}

\let\hat\widehat

\def\Pi{\hat\pi}

\def\lamG{\lam_\Gamma}

\def\mz{m_0}


\Begin{document}


\title{A boundary control problem\\[0.3cm] 
  for the viscous Cahn--Hilliard equation\\[0.3cm]
  with dynamic boundary conditions%
\footnote{{\bf Acknowledgment.}\quad\rm
Some financial support from the
MIUR-PRIN Grant 2010A2TFX2 ``Calculus of variations'' 
and the GNAMPA (Gruppo Nazionale per l'Analisi Matematica, 
la Probabilit\`a e le loro Applicazioni) of 
INdAM (Istituto Nazionale di Alta Matematica) is gratefully 
acknowledged by PC and~GG.}}

\author{}
\date{}
\maketitle
\Bcenter
\vskip-2.5cm
{\large\sc Pierluigi Colli$^{(1)}$}\\
{\normalsize e-mail: {\tt pierluigi.colli@unipv.it}}\\[.25cm]
{\large\sc Gianni Gilardi$^{(1)}$}\\
{\normalsize e-mail: {\tt gianni.gilardi@unipv.it}}\\[.25cm]
{\large\sc J\"urgen Sprekels$^{(2)}$}\\
{\normalsize e-mail: {\tt juergen.sprekels@wias-berlin.de}}\\[.45cm]
$^{(1)}$
{\small Dipartimento di Matematica ``F. Casorati'', Universit\`a di Pavia}\\
{\small via Ferrata 1, 27100 Pavia, Italy}\\[.2cm]
$^{(2)}$
\juerg{
{\small Weierstrass Institute for Applied Analysis and Stochastics}\\
{\small Mohrenstrasse 39, 10117 Berlin, Germany}\\[2mm]
{\small and}\\[2mm]
{\small Department of Mathematics}\\
{\small Humboldt-Universit\"at zu Berlin}\\
{\small Unter den Linden 6, 10099 Berlin, Germany}}\\
[1cm]
\Ecenter

\Begin{abstract}
A boundary control problem for the viscous Cahn--Hilliard equations 
with possibly singular potentials and dynamic boundary conditions
is studied and first order necessary conditions for optimality are proved.
\vskip3mm

\juerg{adjoint state system}.
\End{abstract}

\salta

\pagestyle{myheadings}
\newcommand\testopari{\sc Colli \ --- \ Gilardi \ --- \ Sprekels}
\newcommand\testodispari{\sc Boundary control problem for the viscous Cahn--Hilliard equation}
\markboth{\testodispari}{\testopari}

\finqui


\section{Introduction}
\label{Intro}
\setcounter{equation}{0}

The simplest form of the Cahn--Hilliard equation with or without viscosity (see \cite{CahH, EllSh, EllSt}) reads as follows
\Beq
  \dt y - \Delta w = 0
  \aand
  w = \tau \dt y \gianni{{}-\Delta y +{}} \calW'(y) 
  \quad \hbox{in $\Omega\times (0,T) $}
  \label{cahnh}
\Eeq
where $\Omega$ is the domain where the evolution takes place,
$y$~and $w$ denote the order parameter and the chemical potential, respectively,
and $\tau\geq0$ is the viscosity coefficient. 
Moreover, $\calW'$ represents the derivative of a double well potential~$\calW$,
and typical and important examples are 
the classical regular potential $\calW_{reg}$ 
and the logarithmic double-well potential $\calW_{log}$ given~by
\Bsist
  & \calW_{reg}(r) = \frac14(r^2-1)^2 \,,
  \quad r \in \erre 
  \label{regpot}
  \\
  & \calW_{log}(r) = ((1+r)\ln (1+r)+(1-r)\ln (1-r)) - c r^2 \,,
  \quad r \in (-1,1)
  \label{logpot}
\Esist
where $c>0$ in the latter is large enough in order to kill convexity. 

Moreover, an initial condition like $y(0)=\yz$ 
and suitable boundary conditions must complement the above equations.
As far as the latter are concerned, the most common ones in the literature are
the usual no-flux conditions for both $y$ and~$w$.
However, different boundary conditions have been recently proposed\pier{:
namely,} still the usual no-flux condition for the chemical potential
\Beq
  (\dn w)\suG =0 
  \quad \hbox{on $\Gamma\times (0,T)$}
  \label{nfc}
\Eeq
in order to preserve mass conservation,
and the dynamic boundary condition
\Beq
  (\dn y)\suG + \dt\yG - \DeltaG\yG + \calW'_\Gamma(\yG) = \uG
  \quad \hbox{on $\Gamma\times (0,T)$}
  \label{dbc}
\Eeq
where $\yG$ denotes the trace $y\suG$ on the boundary $\Gamma$ of $\Omega$, 
$\DeltaG$ stands for the Laplace--Beltrami operator on~$\Gamma$, 
$\calW'_\Gamma$ is a nonlinearity analoguous to $\calW'$ 
but now acting on the boundary value of the order parameter, 
and finally $\uG$ is a boundary source term. 
We just quote, among other contributions, \pier{\cite{CFP, Is, MirZelik, PRZ, RZ, WZ}}
and especially the papers \cite{GiMiSchi} and~\cite{CGS}.
In the former, the reader can find the
physical meaning and free energy derivation 
of the boundary value problem given by \eqref{cahnh} and \accorpa{nfc}{dbc},
besides the mathematical \juerg{treatment} of the problem itself.
The latter provides existence, uniqueness and regularity results 
for the same boundary value problem
by assuming that the dominating potential is the boundary potential  $\calW_\Gamma$
rather than the bulk potential~$\calW$
(thus, in contrast to~\cite{GiMiSchi})
and thus \pier{it} is close from this point of view to~\cite{CaCo}, 
where the Allen--Cahn equation with dynamic boundary condition is studied
\pier{(see also~\cite{CF} in which a mass constraint is considered, too)}.

The aim of this paper is \juerg{to study an associated optimal boundary} control problem,
the control being the forcing term $\uG$ that appears on the \rhs\ of the dynamic boundary condition~\eqref{dbc}.
\juerg{While numerous investigations deal with the well-posedness and asymptotic behavior of
Cahn--Hilliard systems, there are comparatively few contributions dedicated to aspects of optimal control.
Usually, these papers treat the non-viscous case ($\tau=0$) and are restricted to distributed controls, with
the no-flux boundary condition $\,(\partial_n y)_{|_\Gamma}=0\,$ assumed in place of the more difficult 
dynamic boundary condition (\ref{dbc}). In this connection, we refer to \cite{WaNa} and \cite{HW1}, where the
latter paper also applies to the case in which the differentiable potentials (\ref{regpot}) or (\ref{logpot}) are 
replaced by the non-differentiable ``double obstacle potential'' given by
}   
\Beq
\juerg{
  {\cal W}_{2obst}(r)= \left\{ 
  \begin{array}{cl}
    c\,(1-r^2) & \hbox{ if } \,|r|\leq 1  
    \\[0.1cm]
    {} + \infty  & \hbox{ if } \,|r| >1  
  \end{array} \right. \quad\,\mbox{for some $c>0$.}
  }
  \label{obstpot}
\Eeq
\juerg{
Note that in this case the second equation in \eqref{cahnh} has to be interpreted as a differential inclusion,
since $\calW'$ cannot be a usual derivative. Instead, the derivative of the convex part of $\calW$ is given by $\partial I_{[-1, 1]}$, 
the subdifferential of the indicator function of the interval $[-1,1]$, which is defined by } 
\Beq
  s \in \partial I_{[-1, 1]} (r) 
  \quad \hbox{ if and only if } \quad 
  s \ \left\{
  \begin{array}{ll}
    \displaystyle
    \leq \, 0 \ &\hbox{if } \ r=-1   
    \\[0.1cm]
    {} = 0 \ &\hbox{if } \ -1< r < 1  
    \\[0.1cm]
    {} \geq 0 \ &\hbox{if } \ r = 1  \,.
    \\[0.1cm]
  \end{array}
  \right. 
  \label{pier5}
\Eeq
\juerg{We remark that the double obstacle case is particularly challenging from the viewpoint of optimal
control, because this case is well known to fall into the class of ``MPEC (Mathematical Programs with Equilibrium
Constraints) Problems''; indeed, the corresponding state system then contains a differential inclusion 
for which all of the standard nonlinear programming constraint qualifications are violated so that
the existence of Lagrange multipliers cannot be shown via standard techniques.} 

\juerg{
Quite recently, also convective Cahn--Hilliard systems have been investigated from the viewpoint of optimal
control. In this connection, we refer to \cite{ZL1} and \cite{ZL2}, where the latter paper deals with the
spatially two-dimensional case. The three-dimensional case with a nonlocal potential \pier{is} studied
in \cite{RoSp}. There also exist contributions dealing with discretized versions of the more general
Cahn--Hilliard--Navier--Stokes system, cf. \cite{HK} and \cite{HW2}. Finally, we
mention the contributions \cite{CGPS} and \cite{CGScon}, in which control problems for
a generalized Cahn--Hilliard system introduced in \cite{Podio} \pier{are} investigated.}

\juerg{To the authors' best knowledge, there are presently no contributions to the optimal boundary control of
viscous or non-viscous Cahn--Hilliard systems with dynamic boundary conditions of the form (\ref{dbc}). \pier{We are aware}, however, \pier{of the} recent contributions \pier{\cite{CS} and \cite{CFS}} for the corresponding Allen--Cahn equation.} 
In particular, \cite{CS}~\pier{treats} both the cases of distributed and boundary controls for logarithmic-type 
potentials as in~\eqref{logpot}.
More precisely, both the existence of optimal controls 
and twice continuous \Frechet\ differentiability for the well-defined control-to-state mapping 
were established, as \juerg{well as first-order necessary 
and second-order sufficient optimality conditions}.
The related paper \cite{CFS} \pier{deals} with the existence of optimal controls 
and the derivation of first-order necessary conditions of optimality
for the more difficult case of the double obstacle potential. The 
method used consists in performing a so-called ``deep quench limit'' 
of the problem studied in~\cite{CS}.

As \juerg{mentioned above, the recent paper} \cite{CGS} \pier{contains} a number of results
that \pier{regard} the problem
obtained by complementing the equations \eqref{cahnh} with
the already \pier{underlined} initial and boundary conditions, namely,
\Bsist
  & \dt y - \Delta w = 0
  \quad \hbox{in $\, Q:=\Omega\times (0,T)$}
  \label{Iprima}
  \\[0.2cm]
  & w = \tau \, \dt y - \Delta y + \calW'(y)
  \quad \hbox{in $\,Q$}
  \label{Iseconda}
  \\[0.2cm]
  & \dn w = 0
  \quad \hbox{on $\, \Sigma:= \Gamma\times (0,T)$}
  \label{Ibc}
  \\[0.2cm]
  & \yG = y\suG
  \aand
  \dt\yG + (\dn y)\suG - \DeltaG\yG + \calW'_\Gamma(\yG) = \uG
  \quad \hbox{on $\, \Sigma$}
  \label{Iterza}
  \\[0.2cm]
  & y(0) = \yz
  \quad \hbox{in $\, \Omega$} .
  \label{Icauchy}
\Esist
\Accorpa\Ipbl Iprima Icauchy
More precisely, existence, uniqueness and regularity results were proved
for general potentials that include \juerg{\accorpa{regpot}{logpot} and~\eqref{obstpot},
and are valid for both the viscous and non-viscous cases}, i.e., by assuming just $\tau\geq0$.
Moreover, further regularity of the solution was ensured 
provided that $\tau>0$ and too singular potentials like \eqref{obstpot} were excluded.
Furthermore, results \juerg{for the linearization around a solution were given as well}.
In such a problem, $\calW'(y)$ and $\calW'_\Gamma(\yG)$ are replaced by $\lam\,y$ and $\lamG\yG$, 
for some given functions $\lam$ and $\lamG$ on $Q$ and~$\Sigma$, respectively.
Therefore, the proper material is already prepared for the control problem \juerg{to be studied here.}

Among several possibilities, \juerg{we choose the tracking-type cost functional }
\Bsist
  & \calJ(y,\yG,\uG)
  & := \frac\bQ 2 \, \norma{y-\zQ}_{L^2(Q)}^2
  + \frac\bS 2 \, \norma{\pier{\yG}-\zS}_{L^2(\Sigma)}^2
  + \frac\bO 2 \, \norma{y(T)-\zO}_{L^2(\Omega)}^2
  \non
  \\
  && \quad +\, \frac\bG 2 \, \norma{\yG(T)-\zG}_{L^2(\Gamma)}^2
  + \frac\bz 2 \, \norma\uG_{L^2(\Sigma)}^2
  \label{Icost}
\Esist
where the functions \juerg{$z_Q, z_\Sigma,z_\Omega,z_\Gamma$}  and the nonnegative constants \juerg{$\bQ,\bS,\bO,\bG,\bz$} are given. \juerg{The control problem then} consists in minimizing \eqref{Icost} 
subject to the state system \Ipbl\ and to the constraint $\uG\in\Uad$, 
where the control box $\Uad$ is given~by
\Bsist
  & \Uad := 
  & \bigl\{ \uG\in\H1\HG\cap\LS\infty:
  \non
  \\
  && \ \uGmin\leq\uG\leq\uGmax\ \aeS,\ \norma{\dt\uG}_2\leq\Mz
  \bigr\}  
  \label{Iuad}
\Esist
for some given functions $\uGmin,\uGmax\in\LS\infty$ and some prescribed positive constant~$\Mz$.

\juerg{In this paper, we confine ourselves to the viscous case $\tau>0$ and avoid potentials like~\eqref{obstpot},
in order to be able to apply all of the results proved in~\cite{CGS}. However,
regular and singular potentials like \eqref{regpot} and \eqref{logpot} are allowed}.
In this framework, we prove both the existence of an optimal control $\uopt$
and first-order necessary conditions for optimality.
To \juerg{this end}, we show the \Frechet\ differentiability of the control-to-state mapping
and introduce and solve a proper adjoint problem, 
which consists in a backward Cauchy problem for the system
\Bsist
  & q = -\Delta p
  \aand
  \gianni{- \dt(p+q) - \Delta q + \lambda q = \phi}
  \quad \hbox{in $Q$}
  \label{Iadj}
  \\
  & \dn p = 0
  \aand
  \dt\qG + \dn q - \DeltaG\qG + \lamG \qG = \phi_\Gamma
  \quad \hbox{on $\Sigma$},
  \label{Ibcadj}
\Esist
where $\qG$ is the trace $q\suG$ of $q$ on the boundary, 
and where the functions $\lambda$, $\lamG$, $\phi$ and $\phi_\Gamma$
are suitably related to the functions \juerg{$z_Q, z_\Sigma,z_\Omega,z_\Gamma$  and the constants $\bQ,\bS,\bO,\bG,\bz$ 
appearing in the cost functional~\eqref{Icost}, as well as} to the state $(\yopt,\yGopt)$ associated to the optimal control~$\uopt$.

The paper is organized as follows.
In the next section, we list our assumptions and state the problem in a precise form.
Moreover, we present some auxiliary material and sketch our results.
The existence of an \pier{optimal} control will be  proved in Section~\ref{OPTIMUM},
while the rest of the paper is devoted 
to the derivation of first-order necessary conditions for optimality.
The final result will be proved in Section~\ref{OPTIMALITY}; it is
prepared in Sections~\ref{FRECHET} and~\ref{ADJOINT},
where we study the control-to-state mapping and solve the adjoint problem.


\section{Statement of the problem and results}
\label{STATEMENT}
\setcounter{equation}{0}

In this section, we describe the problem under study,
present the auxiliary material we need and give an outline of our results.
As in the Introduction,
$\Omega$~is the body where the evolution takes place.
We assume $\Omega\subset\erre^3$
to~be open, bounded, connected, and smooth,
and we write $|\Omega|$ for its Lebesgue measure.
Moreover, $\Gamma$, $\dn$, $\nablaG$ and $\DeltaG$ still stand for
the boundary of~$\Omega$, the outward normal derivative, the surface gradient 
and the Laplace--Beltrami operator, respectively.
Given a finite final time~$T>0$,
we set for convenience
\Bsist
  && Q_t := \Omega \times (0,t)
  \aand
  \Sigma_t := \Gamma \times (0,t)
  \quad \hbox{for every $t\in(0,T]$}
  \label{defQtSt}
  \\
  && Q := Q_T \,,
  \aand
  \Sigma := \Sigma_T \,.
  \label{defQS}
\Esist
\Accorpa\defQeS defQtSt defQS
Now, we \juerg{specify the assumptions on the structure of our system}.
Some of the results we need are known and hold under rather mild conditions.
However, the control problem \juerg{under study} in this paper needs
the high level of regularity for the solution that we are going to specify.
In particular, the values of the state variable have to be bounded
far away from the singularity of the bulk and boundary potentials
in order that the solution to the linearized problem 
\juerg{introduced below} is smooth as well.
Even though all this could be true (for smooth data) also in other situations,
i.e., if the structure of the system is somehow different,
we give a list of assumptions that implies the whole set of conditions listed in~\cite{CGS},
since the latter surely guarantee \juerg{all we need}.
We also assume the potentials to be slightly smoother than in~\cite{CGS},
since this \juerg{will be} useful later on.
In order to avoid a heavy notation,
we write $f$ and $\fG$ in place of $\calW$ and~$\calW_\Gamma$, respectively.
Moreover, as we only consider the case of a positive viscosity coefficient,
we take $\tau=1$ without loss of generality.
\juerg{We make the following assumptions for the structure of
 our system}.
\Bsist
  & -\infty \leq \rmin < 0 < \rmax \leq +\infty
  \label{hpDf}\\[1mm]
    & \hbox{$f,\,\fG:(\rmin,\rmax)\to[0,+\infty)$
    are $C^3$ functions such that}
  \label{hppot}
  \\[1mm]
  & f(0) = \fG(0) = 0
  \aand 
  \hbox{$f''$ and $\fG''$ are bounded from below}
  \label{hpfseconda}
  \\[1mm]
  & \pier{|f'(r)| \leq \eta \,|\fG'(r)| + C
  \quad \hbox{for some $\eta,\, C>0$ and every $r\in(\rmin,\rmax)$}}
  \label{hpcompatib}
  \\[1mm]
  & \pier{\lim\limits_{r\seto\rmin} f'(r)
  = \lim\limits_{r\seto\rmin} \fG'(r) 
  = -\infty 
  \aand
  \lim\limits_{r\neto\rmax} f'(r)
  = \lim\limits_{r\neto\rmax} \fG'(r) 
  = +\infty \, .}
  \label{fmaxmon}
\Esist
\Accorpa\HPstruttura hpDf fmaxmon
We note that \HPstruttura\ imply the possibility
of splitting $f'$ as $f'=\beta+\pi$,
where $\beta$ is a monotone function that diverges at~$r_\pm$
and $\pi$ is a perturbation with a bounded derivative.
Moreover, the same is true for~$\fG$,
so that the assumptions of \cite{CGS} are \juerg{satisfied}.
Furthermore, the choices $f=\calW_{reg}$ and $f=\calW_{log}$
corresponding to \eqref{regpot} and \eqref{logpot} are allowed.

\lastrev{\Brem
\label{referee}
Let us notice that our assumptions \HPstruttura\ perfectly fit the framework 
of the paper \cite{CGS}, in which it is postulated that the boundary potential 
dominates the bulk one via condition \eqref{hpcompatib}. However, the reader 
may wonder whether, especially in cases like $f=\calW_{reg}$, assumption \eqref{hpcompatib} can be released. To this concern, let us refer to the paper \cite{GiMiSchi} where other conditions for the bulk and boundary potentials are 
discussed from the viewpoint of well-posedness: it would be interesting to check 
whether our results for the optimal control problem can be extended. 
\Erem}

Next, in order to \juerg{simplify notation}, we~set
\Bsist
  && V := \Huno, \quad
  H := \Ldue, \quad
  \HG := \LdueG 
  \aand
  \VG := \HunoG 
  \qquad
  \label{defVH}
  \\
  && \calV := \graffe{(v,\vG)\in V\times\VG:\ \vG=v\suG}
  \aand
  \calH := H \times \HG   
  \label{defcalVH}
\Esist
\Accorpa\Defspazi defVH defcalVH
and \juerg{endow these} spaces with their natural norms.
For $1\leq p\leq\infty$, $\norma\cpto_p$ is the usual norm in $L^p$ spaces
($\Lx p$, $\LS p$, etc.),
while $\norma\cpto_X$ stands for the norm in the generic Banach space~$X$.
Furthermore, the symbol $\<\cpto,\cpto>$ stands
for the duality pairing between~$\Vp$, the dual space of~$V$, and $V$ itself.
In the \juerg{following}, it is understood that $H$ is embedded in $\Vp$
in the usual way, i.e.,  \juerg{such that we have
$\<u,v>=(u,v)$ with the inner product $(\,\cdot\,,\,\cdot)$ of~$H$},
for every $u\in H$ and $v\in V$.
Finally, if $u\in\Vp$ and $\underline u\in\L1\Vp$,
we define their \generaliz ed mean values 
$\uO\in\erre$ and $\underline u^\Omega\in L^1(0,T)$ by setting
\Beq
  \uO := \frac 1 {|\Omega|} \, \< u , 1 >
  \aand
  \underline u^\Omega(t) := \bigl( \underline u(t) \bigr)^\Omega
  \quad \aat .
  \label{media}
\Eeq
Clearly, \eqref{media} give the usual mean values when applied to elements
of~$H$ or $\L1H$.

At this point, we can describe the state problem.
For the data we assume~that
\Bsist
  && \lastrev{\yz \in \Hx 2 \,, \quad
  \yz\suG \in \HxG 2}
  \label{hpyz}
  \\
  && \rmin < \yz(x) < \rmax 
  \quad \hbox{for every $x\in\overline\Omega$}.
  \label{hpyzbis}
\Esist
\Accorpa\HPyz hpyz hpyzbis
Moreover, $\uG$ is given in $\H1\HG$.
Even though we could write the equations and the boundary conditions in their strong forms,
we prefer to use the variational formulation of system \Ipbl.
Thus, we look for a triplet $(y,\yG,w)$ \juerg{satisfying}
\Bsist
  & y \in \W{1,\infty}H \cap \H1V \cap \L\infty\Hdue
  \label{regy}
  \\
  & \yG \in \W{1,\infty}\HG \cap \H1\VG \cap \L\infty\HdueG
  \label{regyG}
  \\
  & \yG(t) = y(t)\suG
  \quad \aat
  \label{tracciay}
  \\
  & 
  \rmin < \infess\limits_Q y \leq \supess\limits_Q y < \rmax 
  \label{tuttolinfty}
  \\
  & w \in \L\infty\Hdue
  \label{regw}
\Esist
\Accorpa\Regsoluz regy regw
\juerg{as well as, for almost every $t\in (0,T)$,} the variational equations
\Bsist
  && \iO \pier{\dt} y(t) \, v 
  + \iO \nabla w(t) \cdot \nabla v = 0
  \label{prima}
  \\
  \noalign{\smallskip}
  && \iO w(t) \, \pier{v}
  = \iO \dt y(t) \, v
  + \iG \dt\yG(t) \, \vG
  + \iO \nabla y(t) \cdot \nabla v
  + \iG \nablaG\yG(t) \cdot \nablaG\vG
  \qquad
  \non
  \\
  && \quad {}
  + \iO f'(y(t)) \, v
  + \iG \bigl( \fG'(\yG(t)) - \uG(t) \bigr) \, \vG
  \label{seconda}
\Esist
for every $v\in V$ and every $(v,\vG)\in\calV$, respectively,
and the Cauchy condition
\Beq
  y(0) = \yz \,.
  \label{cauchy}
\Eeq
\Accorpa\State prima cauchy
\Accorpa\Pbl regy cauchy
We note that an equivalent formulation of \accorpa{prima}{seconda}
is given~by
\Bsist
  && \intQ \dt y \, v
  + \intQ \nabla w \cdot \nabla v = 0
  \label{intprima}
  \\
  \noalign{\smallskip}
  && \intQ wv
  = \intQ \dt y \, v
  + \intS \dt\yG \, \vG
  + \intQ \nabla y \cdot \nabla v
  + \intS \nablaG\yG \cdot \nablaG\vG
  \qquad
  \non
  \\
  && \quad {}
  + \intQ f'(y) \, v
  + \intS \bigl( \fG'(\yG) - \uG \bigr) \, \vG
  \qquad
  \label{intseconda}
\Esist
\Accorpa\IntPbl intprima intseconda
for every $v\in\L2 V$ and every $(v,\vG)\in\L2\calV$, respectively.
It is worth \juerg{noting} that (see notation~\eqref{media})
\Bsist
  && (\dt y(t))^\Omega = 0
  \quad \aat
  \aand
  y(t)^\Omega = \mz
  \quad \hbox{for every $t\in[0,T]$}
  \qquad
  \non
  \\
  && \quad \hbox{where $\mz=(\yz)^\Omega$ \juerg{is} the mean value of $\yz$,}
  \label{conserved}
\Esist
as usual for the Cahn--Hilliard equation.

As far as existence, uniqueness, regularity and continuous dependence are concerned,
we can apply Theorems~2.2, 2.3, 2.4, 2.6 and Corollary~2.7 of \cite{CGS}
(where $\calV$ has a slightly different meaning with respect to the present paper)
and~deduce what we need as a particular case. Moreover, as the proofs of the regularity \Regsoluz\ of the solution 
performed in \cite{CGS} mainly rely on a~priori estimates and compactness arguments,
it is clear that a stability estimate holds as well. \lastrev{However, referring to \cite{CGS} let us point out to the reader that the assumption (2.37) explicitly required in the statements of \cite[Thms.~2.4 and 2.6, Cor.~2.7]{CGS} contains the condition $\partial_n {y_{0}}_{{\pier |_\Gamma}}=0$ which is completely useless (actually, it is never employed in the proofs, as one can easily check).}
Therefore, we~have

\Bthm
\label{daCGS}
Assume \HPstruttura\ and \HPyz, and let $\uG\in\H1\HG$.
Then, problem \Pbl\ has a unique solution $(y,\yG,w)$,
and \juerg{this} solution satisfies
\Bsist
  && \norma y_{\W{1,\infty}H \cap \H1V \cap \L\infty\Hdue}
  \non
  \\[1mm]
  && \quad {}
  + \norma\yG_{\W{1,\infty}\HG \cap \H1\VG \cap \L\infty\HdueG}
  \non
  \\[1mm]
  && \quad {}
  + \norma w_{\L\infty\Hdue}
  \,\leq \,c 
  \label{stab}
  \\[2mm]
  && \rmin' \leq y \leq \rmax'
  \quad \aeQ
  \label{faraway}
\Esist
for some constants $c>0$ and $\rmin',\rmax'\in(\rmin,\rmax)$ that depend only on
$\Omega$, $T$, the shape of the nonlinearities $f$ and~$\fG$,
the initial datum~$\yz$,
and on an upper bound for the norm of $\uG$ in $\H1\HG$. 
Moreover, if $u_{\Gamma\!,i}\in\H1\HG$, $i=1,2$, are two forcing terms
and $(y_i,y_{\Gamma\!,i},w_i)$ are the corresponding solutions,
the inequality
\Bsist
  && \norma{y_1-y_2}_{\L\infty H}^2
  + \norma{y_{\Gamma\!,1}-y_{\Gamma\!,2}}_{\L\infty\HG}^2
  \non
  \\[1mm]
  && \quad {}
  + \norma{\nabla(y_1-y_2)}_{\L2H}^2
  + \norma{\nablaG(y_{\Gamma\!,1}-y_{\Gamma\!,2)}}_{\L2\HG}^2
  \non
  \\[1mm]
  && \leq c \, \norma{u_{\Gamma\!,1}-u_{\Gamma\!,2}}_{\L2\HG}^2
  \label{contdip}
\Esist
holds true with a constant $c$ that depends only on
$\Omega$, $T$, and the shape of the nonlinearities $f$ and~$\fG$.
\Ethm

Once well-posedness for problem \Pbl\ is established,
we can address the corresponding control problem.
As in the Introduction, given four functions
\Beq
  \zQ \in \LQ2 , \quad
  \zS \in \LS2 , \quad
  \zO \in \Lx2
  \aand
  \zG \in \LxG 2
  \label{hpzzzz}
\Eeq
and nonnegative constants $\bQ,\, \bS,\, \bO,\, \bG,\, \bz$,
we~set
\Bsist
  & \calJ(y,\pier{y_\Gamma,} \uG)
  & := \frac\bQ 2 \, \norma{y-\zQ}_{L^2(Q)}^2
  + \frac\bS 2 \, \norma{\juerg{y_\Gamma}-\zS}_{L^2(\Sigma)}^2
  + \frac\bO 2 \, \norma{y(T)-\zO}_{L^2(\Omega)}^2
  \non
  \\
  && {} + \frac\bG 2 \, \norma{\yG(T)-\zG}_{L^2(\Gamma)}^2
  + \frac\bz 2 \, \norma\uG_{L^2(\Sigma)}^2
  \label{defcost}
\Esist
for, say, $y\in\C0H$, $\yG\in\C0\HG$ and $\uG\in\LS2$,
and consider the problem of minimizing the cost functional \eqref{defcost}
subject to the constraint $\uG\in\Uad$,
where the control box $\Uad$ is given~by
\Bsist
  & \Uad := 
  & \bigl\{ \uG\in\H1\HG\cap\LS\infty:
  \non
  \\
  && \ \uGmin\leq\uG\leq\uGmax\ \aeS,\ \norma{\dt\uG}_2\leq\Mz
  \bigr\} 
  \label{defUad}
\Esist
and to the state system \State.
\gianni{We simply assume~that
\Beq
  \Mz > 0 , \quad
  \uGmin ,\, \uGmax \in \LS\infty 
  \aand
  \hbox{$\Uad$ is nonempty}.
  \label{hpUad}
\Eeq
}%
Here is our first result.

\Bthm
\label{Optimum}
Assume \HPstruttura\ and \HPyz,
and let $\Uad$ and $\calJ$ be defined by \eqref{defUad} and \eqref{defcost}
under the \juerg{assumption}s \eqref{hpUad} and~\eqref{hpzzzz}.
Then, there exists $\uopt\in\Uad$ such~that
\Beq
  \calJ(\yopt,\yGopt,\uopt)
  \leq \calJ(y,\yG,\uG)
  \quad \hbox{for every $\uG\in\Uad$}
  \label{optimum}
\Eeq
where $\yopt$, $\yGopt$, $y$ and $\yG$
are the components of the solutions $(\yopt,\yGopt,\wopt)$ and $(y,\yG,w)$
to~the state system \Pbl\ corresponding to the controls
$\uopt$ and~$\uG$, respectively.
\Ethm

From now on, it is understood that the assumptions \HPstruttura\ and \HPyz\
on the structure and on the initial datum $\yz$ are satisfied
and that the cost functional and the control box are defined by \eqref{defcost} and~ \eqref{defUad}
under the \juerg{assumptions} \eqref{hpUad} and~\eqref{hpzzzz}.
Thus, we do not remind anything of that in the statements given in the sequel.

\juerg{Our next aim is to formulate} necessary optimality conditions.
To this end, by recalling the \pier{involved definitions} {\Defspazi},
we introduce the control-to-state mapping by 
\Bsist
  && \calX := \H1\HG \cap \LS\infty
  \aand
  \calY := \H1\calH \cap \L\infty\calV
  \label{defXY}
  \\
  && \hbox{$\calU$ is an open set in $\calX$ that includes $\Uad$}
  \label{defU}
  \\
  && \calS : \calU \to \calY, \quad
  \uG \mapsto \calS(\uG) =: (y,\yG),
  \quad \hbox{where }\,\hbox{$(y,\yG,w)$ is the unique}
  \non
  \\
  &&  \mbox{solution to \Pbl\ corresponding to }\,\uG\,,
  \qquad
  \label{defS}
\Esist
\Accorpa\Controltostate defXY defS
as well as the so-called ``reduced cost functional''
\Beq
  \redJ : \calU \to \erre,
  \quad \hbox{defined by} \quad
  \redJ(\uG) := \calJ(y,\yG,\uG)
  \quad \hbox{where} \quad
  (y,\yG) = \calS(\uG).
  \label{defredJ}
\Eeq
As $\Uad$ is convex, the desired necessary condition for optimality~is
\Beq
  \< D\redJ(\uopt) , \vG-\uopt > \geq 0
  \quad \hbox{for every $\vG\in\Uad$}
  \label{precnopt}
\Eeq
provided that the derivative $D\redJ(\uopt)$ exists in the dual space $(\H1\HG)^*$
at least in the G\^ateaux sense.
Then, the natural approach consists in proving that 
$\calS$ is \Frechet\ differentiable at $\uopt$
and applying the chain rule.
As we \juerg{shall} see in Section~\ref{FRECHET}, 
this leads to the linearized problem that we describe at once
and that can be stated starting from a generic element $\uG\in\calU$.
Let $\uG\in\calU$ and $\hG\in\H1\HG$ be given.
We set $(y,\yG):=\calS(\uG)$ and
\Beq
  \lambda := f''(y)
  \aand
  \lamG := \fG''(\yG) .
  \label{deflam}
\Eeq
Then the problem consists in finding $(\xi,\xiG,\eta)$ 
satisfying the analogue of the regularity requirements \Regsoluz,  
solving \aat\ the variational equations
\Bsist
  && \iO \dt\xi(t) \, v 
  + \iO \nabla\eta(t) \cdot \nabla v = 0
  \label{linprima}
  \\
  \noalign{\smallskip}
  && \iO \eta(t) v
  = \iO \dt\xi(t) \, v
  + \iG \dt\xiG(t) \, \vG
  + \iO \nabla\xi(t) \cdot \nabla v
  + \iG \nablaG\xiG(t) \cdot \nablaG\vG
  \qquad
  \non
  \\
  && \quad {}
  + \iO \lambda(t) \, \xi(t) \, v
  + \iG \bigl( \lamG(t) \, \xiG(t) - \hG(t) \bigr) \, \vG
  \label{linseconda}
\Esist
for every $v\in V$ and every $(v,\vG)\in\calV$, respectively,
and \gianni{satisfying} the Cauchy condition
\Beq
  \xi(0) = 0 \,.
  \label{lincauchy}
\Eeq
\Accorpa\Linpbl linprima lincauchy
Note that property \eqref{conserved} applied to $\xi$ becomes
\Beq
  \xi^\Omega(t) = 0
  \quad \aat,
  \label{xiconserved}
\Eeq
since $\xi(0)=0$.
As $\lambda$ and $\lamG$ are bounded,
we can apply \cite[Cor.~2.5]{CGS} \juerg{to conclude the following result.}

\Bprop
\label{Existlin}
Let $\uG\in\calU$.
Moreover, let $(y,\yG)=\calS(\uG)$ and define $\lambda$ and $\lamG$ by~\eqref{deflam}.
Then, for every $\hG\in\H1\HG$,
there exists a unique triplet $(\xi,\xiG,\eta)$
satisfying the analogue of \Regsoluz\
and solving the linearized problem \Linpbl.
\Eprop

\gianni{Namely}, we shall prove that the \Frechet\ derivative 
$D\calS(\uG)\in\calL(\calX,\calY)$ 
actually exists \pier{and the value that it} \juerg{assigns} to the generic element $\hG\in\calX$
is precisely $(\xi,\xiG)\in\calY$, where $(\xi,\xiG,\eta)$ is
the solution to the linearized problem corresponding to the datum~$\hG$.
This will be done in Section~\ref{FRECHET}.
Once this \juerg{will be established, we may use the chain rule to prove}
that the necessary condition \eqref{precnopt} for optimality takes the form
\Bsist
  && \bQ \intQ (\yopt - \zQ) \xi
  + \bS \intS (\yGopt - \zS) \xiG
  + \bO \iO (\yopt(T) - \zO) \xi(T)
  \non
  \\
  && \quad {}
  + \bG \iG (\yGopt(T) - \zG) \xiG(T)
  + \bz \intS \uopt (\vG - \uopt)
  \geq 0
  \quad \hbox{for any $\vG\in\Uad$},
  \qquad
  \label{cnopt}
\Esist
where, for any given $\vG\in\Uad$, 
the functions $\xi$ and $\xiG$ are the first two components
of the solution $(\xi,\xiG,\eta)$ to the linearized problem
corresponding to $\hG=\vG-\uopt$.

The final step then consists in eliminating the pair $(\xi,\xiG)$ from~\eqref{cnopt}.
This will be done by introducing a triplet $(p,q,\qG)$ 
that fulfills the regularity requirements
\Bsist
  && p \in \H1{\Hx2} \cap \L2{\Hx4}
  \label{regp}
  \\
  && q \in \H1H \cap \L2{\Hx2}
  \label{regq}
  \\
  && \qG \in \H1\HG \cap \L2{\HxG2}
  \label{regqG}
  \\
  && \qG(t) = q(t)\suG
  \quad \aat
  \label{qqG}
\Esist
\Accorpa\Regadj regp qqG
and solves \pier{a suitable} \juerg{backward-in-time problem (the so-called ``adjoint system'')}: \pier{namely,}
the variational equations
  \Beq
  \iO q(t) \, v
  = \iO \nabla p(t) \cdot \nabla v
  \label{primaadj} \pier{\quad \hbox{for all } v\in V \hbox{ and } t\in [0,T]} 
  \Eeq
  \begin{align}
  \gianni{- \iO \dt \bigl( p(t) + q(t) \bigr) \, v
  + \iO \nabla q(t) \cdot \nabla v
  + \iO f''(\yopt(t)) \, q(t) \, v}
  \qquad
  \non
  \\
  \quad {}
  - \iG \dt\qG(t) \, \vG
  + \iG \nablaG\qG(t) \cdot \nablaG\vG
  + \iG \fG''(\yGopt(t)) \, \qG(t) \, \vG
  \non
  \\
  = \iO \bQ \bigl( \yopt(t) - \zQ(t) \bigr) v
  + \iG \bS \bigl( \yGopt(t) - \zS(t) \bigr) \vG
  \non \\
  \pier{\quad \hbox{for every } (v,\vG)\in\calV \hbox{ and a.a. } t\in (0,T)}
  \label{secondaadj}
\end{align}
\gianni{and the final condition}
\begin{align}
  \iO (p+q)(T) \, v 
  + \iG \qG(T) \, \vG
  = \iO \bO \bigl( \yopt(T) - \zO \bigr) v
  + \iG \bG \bigl( \yGopt(T) - \zG \bigr) \vG \non \\ 
  \pier{\quad \hbox{for every } (v,\vG)\in\calV}  \label{cauchyadj} 
\end{align}
\Accorpa\Pbladj primaadj cauchyadj
\pier{have to be} satisfied.
\gianni{\pier{Some} assumptions will be given in order that
this problem has a unique solution, and the optimality condition~\eqref{cnopt}
will be rewritten in a much simpler form.
For instance, one can assume that
\Beq
  \bO = 0 \aand \bG = 0 \,.
  \label{hpzz}
\Eeq
In Sections~\ref{ADJOINT} and~\ref{OPTIMALITY}, 
we will prove the results stated below and sketch how to avoid \eqref{hpzz} 
by weakening a little the summability requirements on the solution
(see the forthcoming Remark~\ref{Pesi}).}

\Bthm
\label{Existenceadj}
\gianni{Let $\uopt$ and $(\yopt,\yGopt)=\calS(\uopt)$
be an optimal control and the corresponding state
and assume in addition that \eqref{hpzz} holds}.
Then the adjoint problem \Pbladj\ has a unique solution $(p,q,\qG)$
satisfying the regularity conditions \Regadj.
\Ethm

\Bthm
\label{CNoptadj}
Let $\uopt$ be an optimal control.
Moreover, let $(\yopt,\yGopt)=\calS(\uopt)$ and $(p,q,\qG)$
be the associate state and the unique solution to the adjoint problem~\Pbladj\
given by Theorem~\ref{Existenceadj}.
Then we~have
\Beq
  \intS (\qG + \bz \uopt) (\vG - \uopt) \geq 0
  \quad \hbox{for every $\vG\in\Uad$}.
  \label{cnoptadj}
\Eeq
\Ethm

\gianni{In particular, if $\bz>0$, \pier{we remark that} $\uG$ is just a projection,
namely
\Beq
  \hbox{\sl\mathsurround 3pt
  $\uopt$ is the orthogonal projection of $-\qG/\bz$ on $\Uad$}
  \label{projection}
\Eeq
with respect to the standard scalar product in~$\LS2$.}

\juerg{In the remainder of this section, we recall some well-known facts
and introduce some notation.}
First of all, we often owe to the elementary Young inequality
\Beq
  ab \leq \delta a^2 + \frac 1 {4\delta} \, b^2
  \quad \hbox{for every $a,b\geq 0$ and $\delta>0$}
  \label{young}
\Eeq
and to the \Holder\ inequality.
Moreover, we account for the well-known Poincar\'e inequality
\Beq
  \normaV v^2 \leq C \bigl( \normaH{\nabla v}^2 + |\vO|^2 \bigr)
  \quad \hbox{for every $v\in V$}
  \label{poincare}
\Eeq
where $C$~depends only on~$\Omega$.
Next, we recall a tool that is generally used 
in the context of problems related to the Cahn--Hilliard equations.
We define
\Beq
  \dom\calN := \graffe{\vstar\in\Vp: \ \vstar^\Omega = 0}
  \aand
  \calN : \dom\calN \to \graffe{v \in V : \ \vO = 0}
  \label{predefN}
\Eeq
by setting for $\vstar\in\dom\calN$
\Beq
  {\calN\vstar \in V, \quad
  (\calN\vstar)^\Omega = 0 ,
  \aand
  \iO \nabla\calN\vstar \cdot \nabla z = \< \vstar , z >
  \quad \hbox{for every $z\in V$}}
  \label{defN}
\Eeq
i.e., $\calN\vstar$ is the solution $v$ to the \generaliz ed Neumann problem for $-\Delta$
with datum~$\vstar$ that satisfies~$\vO=0$.
Indeed, if $\vstar\in H$, the above variational equation means
$-\Delta\calN\vstar = \vstar$ and $\dn\calN\vstar = 0$.
As $\Omega$ is bounded, smooth, and connected,
it turns out that \eqref{defN} yields a well-defined isomorphism
which also satisfies
\Bsist
  && \calN\vstar\in\Hx{s+2}
  \aand
  \norma{\calN\vstar}_{\Hx{s+2}}
  \leq C_s \norma\vstar_{\Hx s}
  \non
  \\
  && \quad \hbox{if $s\geq0$}
  \aand
  \vstar \in \Hx s \cap \dom\calN 
  \label{regN}
\Esist
with a constant $C_s$ that depends only on $\Omega$ and~$s$.
Moreover, we have
\Beq
  \< \ustar , \calN \vstar >
  = \< \vstar , \calN \ustar >
  = \iO (\nabla\calN\ustar) \cdot (\nabla\calN\vstar)
  \quad \hbox{for $\ustar,\vstar\in\dom\calN$}
  \label{simmN}
\Eeq
whence also
\Beq
  2 \< \dt\vstar(t) , \calN\vstar(t) >
  = \frac d{dt} \iO |\nabla\calN\vstar(t)|^2
  = \frac d{dt} \, \normaVp{\vstar(t)}^2
  \quad \aat
  \label{dtcalN}
\Eeq
for every $\vstar\in\H1\Vp$ satisfying $(\vstar)^\Omega=0$ \aet.

We conclude this section by stating a general rule
we use as far as constants are concerned,
in order to avoid a boring notation.
Throughout the paper,
the small-case symbol $c$ stands for different constants which depend only
on~$\Omega$, on the final time~$T$, the shape of the nonlinearities
and on the constants and the norms of
the functions involved in the assumptions of our statements.
Hence, the meaning of $c$ might
change from line to line and even in the same chain of equalities or inequalities.
On the contrary, capital letters (with or without subscripts)
stand for precise constants which we can refer~to.


\section{Existence of an optimal control}
\label{OPTIMUM}
\setcounter{equation}{0}

We prove Theorem~\ref{Optimum} by the direct method,
\juerg{recalling  that $\Uad$ is a nonempty, closed and} convex set in~$\LS2$.
Let $\graffe{\unG}$ be a minimizing sequence for the optimization problem
and, for any~$n$, let us take the corresponding solution $(\yn,\ynG,\wn)$ to problem~\Pbl.
Thus, $\unG\in\Uad$ for every~$n$,
and we can account for the definition \eqref{defUad} of~$\Uad$,
assumptions~\eqref{hpUad}
and \pier{estimates~\eqref{stab}--\eqref{faraway}} for the solutions. 
In particular, we have~that
\Beq
  \rmin < \rmin' \leq \yn \leq \rmax' < \rmax
  \quad \hbox{\aeQ\ and for every $n$}.
  \label{dafaraway}
\Eeq
Next, \pier{owing to} weak star and strong compactness results
(see, e.g., \cite[Sect.~8, Cor.~4]{Simon}),
we~deduce that suitable (not relabeled) subsequences exist such~that
\Bsist
  & \unG \to \uopt
  & \quad \hbox{weakly star in $\LS\infty\cap\H1H$}
  \non
  \\
  & \yn \to \yopt
  & \quad \hbox{weakly star in $\W{1,\infty}H \cap \H1V \cap \L\infty\Hdue$}
  \non
  \\
  && \qquad \hbox{and strongly in $\C0V$}
  \non
  \\
  & \ynG \to \yGopt
  & \quad \hbox{weakly star in $\W{1,\infty}\HG \cap \H1\VG \cap \L\infty\HdueG$}
  \non
  \\
  && \qquad \hbox{and strongly in $\C0\VG$}
  \non
  \\
  & \wn \to \wopt
  & \quad \hbox{weakly star in $\L\infty\Hdue$}.
  \non
\Esist
Clearly, $\uopt\in\Uad$.
Moreover, $\yopt(0)=\yz$.
Furthermore, by \eqref{dafaraway} and the regularity of $f$ and $\fG$ 
we have assumed in \HPstruttura,
we also deduce that $f'(\yn)$ and $\fG'(\ynG)$
converge to $f'(\yopt)$ and $\fG'(\yGopt)$, e.g., strongly in $\C0H$ and $\C0\HG$, respectively.
Thus, we can pass to the limit in the integrated variational formulation~\IntPbl\
written for $(\yn,\ynG,\wn)$ and $\unG$ 
and immediately conclude that the triplet $(\yopt,\yGopt,\wopt)$ 
is the solution $(y,\yG,w)$ to \IntPbl\ corresponding to~$\uG:=\uopt$.
Finally, by semicontinuity,
it is clear that the value $\calJ(\yopt,\yGopt,\uopt)$ 
is the infimum of the cost functional since we have started from a minimizing sequence.
\qed


\section{The control-to-state mapping}
\label{FRECHET}
\setcounter{equation}{0}

We recall the definitions \Controltostate\ of the spaces \pier{$\calX$, $\calY$, the set
$\calU$ and} the map~$\calS$.
As \pier{sketched} in Sections~\ref{STATEMENT}, 
the main point is the \Frechet\ differentiability of the control-to-state mapping~$\calS$.
Our result on that point is prepared by a stability estimate
given by the following lemma.

\Blem
\label{Stability}
Let $u_{\Gamma,i}\in\H1\HG$ for $i=1,2$
and let $(y_i,y_{\Gamma,i},w_i)$ be the corresponding solutions
given by Theorem~\ref{daCGS}.
Then, the following estimate
\Beq
  \norma{(y_1,y_{\Gamma,1}) - (y_2,y_{\Gamma,2})}_{\calY}
  \leq c \norma{u_{\Gamma,1} - u_{\Gamma,2}}_{\L2\HG}
  \label{stability}
\Eeq
holds true for some constant $c>0$ that depends only on
$\Omega$, $T$, the shape of the nonlinearities $f$ and~$\fG$,
and the initial datum~$\yz$.
\Elem

\Bdim
We set for convenience
\Beq
  \uG := u_{\Gamma,1} - u_{\Gamma,2} \,, \quad
  y := y_1 - y_2 \,, \quad
  \yG := y_{\Gamma,1} - y_{\Gamma,2}
  \aand
  w := w_1 - w_2 \,.
  \non
\Eeq
By writing problem \Pbl\ for both solutions $(y_i,y_{\Gamma,i},w_i)$
and taking the difference, we immediately derive~that
\Bsist
  && \iO \dt y(t) \, v 
  + \iO \nabla w(t) \cdot \nabla v = 0
  \label{primadiff}
  \\
  \noalign{\smallskip}
  && \iO w(t) \, v
  = \iO \dt y(t) \, v
  + \iG \dt\yG(t) \, \vG
  + \iO \nabla y(t) \cdot \nabla v
  + \iG \nablaG\yG(t) \cdot \nablaG\vG
  \qquad
  \non
  \\
  && \quad {}
  + \iO \bigl( f'(y_1(t)) - f'(y_2(t)) \bigr) \, v
  + \iG \bigl( \fG'(y_{\Gamma,1}(t)) - \fG'(y_{\Gamma,2}(t)) - \uG(t) \bigr) \, \vG
  \label{secondadiff}
\Esist
\aat\ and for every $v\in V$ and every $(v,\vG)\in\calV$, respectively.
Moreover, $y(0)=0$ and $\dt y$ has zero mean value 
since \eqref{conserved} holds for~$\dt y_i$.
Therefore, $\calN\dt y$ is well defined \aet\ (see~\eqref{predefN})
and we can test \accorpa{primadiff}{secondadiff} written at the time~$s$
by $\calN(\dt y(s))$ and $-\dt(y,\yG)(s)$, respectively.
Then we add the \juerg{resulting equalities}
and integrate over $(0,t)$ with respect to~$s$,
where $t\in(0,T)$ is arbitrary.
We obtain
\Bsist
  && \intQt \dt y \, \juerg{\calN(\dt y)}
  + \intQt \juerg{\nabla w }\cdot \nabla\calN(\dt y)
  - \intQt w \, \dt y
  \non
  \\
  && \quad {}
  + \intQt |\dt y|^2
  + \intSt |\dt\yG|^2
  + \frac 12 \iO |\nabla y(t)|^2
  + \frac 12 \iG |\nablaG\yG(t)|^2
  \non
  \\
  && = - \intQt \bigl( f'(y_1) - f'(y_2) \bigr) \dt y
  - \intSt \bigl( \fG'(y_{\Gamma,1}) - \fG'(y_{\Gamma,2}) \bigr) \dt\yG 
  + \intSt \uG \, \dt\yG \,.
  \qquad
  \label{perstab}
\Esist
By accounting for \juerg{\eqref{simmN}} and \eqref{defN}, we have
\Beq
  \intQt \dt y \, \calN(\dt y)
  + \intQt \juerg{\nabla w} \cdot \nabla\calN(\dt y)
  - \intQt w \, \dt y
  = \juerg{ \intQt |\nabla\calN\dt y|^2} \,  \geq \,0 \,.
  \non
\Eeq
Moreover, all the other integrals on the \lhs\ of~\eqref{perstab} are nonnegative.
The first two terms on the \rhs\ need the same treatement
and we only deal with the first of them.
We notice that both $y_1$ and $y_2$ satisfy~\eqref{faraway}
and that \juerg{$f'$} is \Lip\ continuous on~$[\rmin',\rmax']$.
By using this and the \Holder, Young and Poincar\'e inequalities 
(see~\eqref{poincare}), we derive~that
\Beq
  - \intQt \bigl( f'(y_1) - f'(y_2) \bigr) \dt y
  \leq \frac 14 \intQt |\dt y|^2
  + c \intQt |\nabla y|^2 \,.
  \non
\Eeq
Finally, we simply have
\Beq
  \intSt \uG \, \dt\yG
  \leq \frac 14 \intSt |\dt\yG|^2
  + \intSt |\uG|^2 \,.
  \non
\Eeq
By combining these inequalities and \eqref{perstab},
we see that we can apply the standard Gronwall lemma.
This directly yields~\eqref{stability}.
\Edim

\Bthm
\label{Fdiff}
Let $\uG\in\calU$.
Moreover, let $(y,\yG)=\calS(\uG)$ and define $\lambda$ and $\lamG$ by~\eqref{deflam}.
Then \juerg{the control-to-state mapping $\calS:{\cal U}\subset {\cal X}\to{\cal Y} $  is \Frechet\ differentiable at $\uG$, 
and its \Frechet\ derivative $D\calS(\uG) \in {\cal L}({\cal X},{\cal Y})$ is given as follows:
for $\hG\in\calX$, the value of $D\calS(\uG)$ at $\hG$ is the pair $(\xi,\xiG)$,
where $(\xi,\xiG,\eta)$ is the unique solution to the linearized problem \Linpbl.}
\Ethm

\Bdim
\juerg{At first, a closer inspection of the proof of Theorem 4.1 in \cite{CGS} for the linear case reveals 
that the linear mapping, which assigns to each $h_\Gamma\in {\cal X}$ the pair $(\xi,\xi_\Gamma)$, where 
$(\xi,\xi_\Gamma,\eta)$ is the associated unique solution  to the linearized system \Linpbl, is bounded as 
a mapping from ${\cal X}$ into ${\cal Y}$. Hence, if $D{\cal S}(u_\Gamma)$ has in fact the asserted form, then
it belongs to ${\cal L}({\cal X},{\cal Y})$.}

In the following, it understood that $\norma\hG_{\calX}$
is small enough in order that $\uG+\hG\in\calU$.
As we would like writing the inequality that shows the desired differentiability in a simple form,
we introduce some auxiliary functions.
First of all, we also need the third component $w$ 
of the solution $(y,\yG,w)$ associated to~$\uG$.
Moreover, given~$\hG\in\calX$ small enough, we~set
\Bsist
  && (\yh,\yhG,\wh) := \hbox{solution to \Pbl\ corresponding to $\uG+\hG$},
  \non
  \\
  && \quad \hbox{whence} \quad
  (\yh,\yhG) = \calS(\uG+\hG) 
  \non
  \\
  && \qh := \yh - y - \xi , \quad
  \qhG := \yhG - \yG - \xiG 
  \aand
  \zh := \wh - w - \eta \,.
  \non
\Esist
By \juerg{the definition of the notion  \Frechet\ derivative, we need to show that
$\norma{(\qh,\qhG)}_{\calY}=o(\norma\hG_{\calX})$
as $\norma\hG_{\calX}\to0$.
We prove a preciser} estimate, namely
\Beq
  \norma{(\qh,\qhG)}_{\calY}
  \leq c \norma\hG_{\LS2}^2 \,.
  \label{tesiFrechet}
\Eeq
By definition, the triplets \pier{$(\yh,\yhG,\wh)$ and $(y,\yG,w)$}
satisfy problem \Pbl\ with data \pier{$\uG+ \hG$ and $\uG$}, respectively.
Moreover, $(\xi,\xiG,\eta)$ solves the linearized problem~\Linpbl.
By writing everything and taking the difference, we obtain~\aat
\Bsist
  &&\iO \dt\qh(t) \, v 
  + \iO \nabla \zh(t) \cdot \nabla v = 0
  \label{primah}
  \\[0.2cm]  
  && \iO \zh(t) \, \pier{v}
  = \iO \dt\qh(t) \, v
  + \iG \dt\qhG(t) \, \vG
  + \iO \nabla\qh(t) \cdot \nabla v
  + \iG \nablaG\qhG(t) \cdot \nablaG\vG
  \qquad
  \non
  \\
  && \quad {}
  + \iO \bigl( f'(\yh(t)) - f'(y(t)) - f''(y(t)) \xi(t) \bigr) \, v
  \qquad
  \non
  \\
  && \quad {}
  + \iG \bigl( \fG'(\yhG(t)) - \fG'(\yhG(t)) - \fG''(\yG(t)) \xiG(t) \bigr) \, \vG 
  \label{secondah}
\Esist
for every $v\in V$ and every $(v,\vG)\in\calV$, respectively.
Moreover, $\qh(0)=0$.
We observe that the choice $v=1$ in \eqref{primah}
yields that $\dt\qh(t)$ has zero mean value and thus belongs to the domain of~$\calN$ \aat\
(see \eqref{predefN}).
Therefore, we can test \accorpa{primah}{secondah} written at the time~$s$
by $\calN(\dt\qh(s))$ and $-\dt(\qh,\qhG)(s)$, respectively.
Then, we add the \juerg{resulting equalities}
and integrate over $(0,t)$ with respect to~$s$,
where $t\in(0,T)$ is arbitrary.
We obtain
\Bsist
  && \intQt \dt\qh \, \calN(\dt\qh)
  + \intQt \nabla\zh \cdot \nabla\calN(\dt\qh)
  - \intQt \zh \, \dt\qh
  \non
  \\
  && \quad {}
  + \intQt |\dt\qh|^2
  + \intSt |\dt\qhG|^2
  + \frac 12 \iO |\nabla\qh(t)|^2
  + \frac 12 \iG |\nablaG\qhG(t)|^2
  \non
  \\
  && = - \intQt \bigl( f'(\yh) - f'(y) - f''(y) \xi \bigr) \, \dt\qh
  - \intSt \bigl( \fG'(\yhG) - \fG'(\yhG) - \fG''(\yG) \xiG \bigr) \, \dt\qhG \,.
  \qquad
  \label{perFrechet}
\Esist
As in the proof of Lemma~\ref{Stability},
the sum of the first three integrals on the \lhs\ of~\eqref{perFrechet} is nonnegative
as well as each of the other terms.
Now, we estimate the first integral on the \rhs. 
We write the second order Taylor expansion of the $C^2$ function~$f'$ (see~\eqref{hppot})  
at~$y$ in the Lagrange form.
As $\yh-y=\xi+\qh$, we obtain
\Beq
  f'(\yh) - f'(y) - f''(y) \xi
  = f''(y) \, \qh + \frac 12 \, f'''(\sigma) |\yh-y|^2,
  \non 
\Eeq
with some function $\sigma$ taking its values between the ones of $\yh$ and~$y$.
As $\yh$ and $y$ are bounded away from~$r^\pm$ 
(see~\eqref{faraway}, which holds for both of them),
$f''(y)$ and $f'''(\sigma)$ are \juerg{bounded} in $\LQ\infty$, 
and the above expansion yields
\Beq
  |f'(\yh) - f'(y) - f''(y) \xi|
  \leq c \bigl( |\qh| + |\yh-y|^2 \bigr) .
  \non
\Eeq
Hence, we have
\Beq
  - \intQt \bigl( f'(\yh) - f'(y) - f''(y) \xi \bigr) \, \dt\qh
  \leq C_1 \intQt |\qh| \, |\dt\qh|
  + C_2 \intQt |\yh-y|^2 |\dt\qh| 
  \label{dataylor}
\Eeq
where we have marked the constants in front of the last two integrals
for a future reference.
We deal with the first term on the \rhs\ of the last inequality as follows:
\Beq
  C_1 \intQt |\qh| \, |\dt\qh|
  \leq \frac 14 \intQt |\dt\qh|^2
  + c \intQt |\qh|^2  
  \leq \frac 14 \intQt |\dt\qh|^2
  + c \intQt |\nabla\qh|^2  \,,
  \non
\Eeq
by the Poincar\'e inequality \eqref{poincare},
since $\qh\pier{{}= \yh - y - \xi}$ has zero mean value.
Indeed, $(\yh)^\Omega=y^\Omega$ and $\xi^\Omega=0$ 
since $\yh(0)=y(0)$ and $\xi(0)=0$
(see~\eqref{conserved}).
As far as the last term in \eqref{dataylor} is concerned, we can estimate it this way
\Bsist
  && C_2 \intQt |\yh-y|^2 |\dt\qh|
  \leq c \norma{\yh-y}_{\L\infty V}^2 \Bigl( \intQt |\dt\qh|^2 \Bigr)^{1/2}
  \non
  \\
  && \leq \frac 14 \intQt |\dt\qh|^2
  + c \norma{\yh-y}_{\L\infty V}^4 
  \leq \frac 14 \intQt |\dt\qh|^2
  + c \norma\hG_{\LS2}^4
  \non
\Esist
thanks to the stability estimate~\eqref{stability}.
As the same calculation can be done 
for the last term on the \rhs\ of~\eqref{perFrechet},
we can combine, apply the standard Gronwall lemma,
and conclude that \eqref{tesiFrechet} holds true.

\Edim


\section{The adjoint problem}
\label{ADJOINT}
\setcounter{equation}{0}

In this section, we prove Theorem~\ref{Existenceadj}, i.e.,
we show that problem \Pbladj\ has a unique solution under the further assumptions~\eqref{hpzz}.
\gianni{Moreover, we briefly show how \eqref{hpzz} can be avoided by just requiring less regularity to the solution
(see Remark~\ref{Pesi})}.

In order to solve problem \Pbladj,
we first prove that it is equivalent to a decoupled problem
that can be solved by first finding $q$ and then \gianni{reconstructing~$p$}.
The basic ideas are explained at once.
We note that the function $q(t)$ has zero mean value \aat,
as we immediately see by choosing $v=1$ in~\eqref{primaadj}.
So, if we introduce the mean value function $\pO\in C^0([0,T])$
(see~\eqref{media}),
we~realize that, \aat, $(p-\pO)(t)$ satisfies definition~\eqref{defN}
with $\vstar=q(t)$. We thus have
\gianni{$p(t)-\pO(t)=\calN(q(t))$}.
On the other hand, for any fixed~$t$, the function $\pO(t)$ is a constant;
thus, it is orthogonal in $\Lx2$
to the subspace of functions having zero mean value.
\gianni{Thus, $p$~is completely eliminated} from equation \eqref{secondaadj}
if we confine ourselves to use test functions with zero mean value.
Similar remarks have to be done for the final condition on $p+q$
that appears in~\eqref{cauchyadj}.
Whenever we find a solution $(q,\qG)$ to this new problem, then
\gianni{we can reconstruct $p$ as just said,
provided that we can calculate~$\pO$}.
All this is made precise in our next theorem.
As we are going to use test functions with zero mean value,
we introduce the proper spaces by
\Beq
  \HO := \graffe{(v,\vG) \in \calH :\ \vO=0}
  \aand
  \VO := \HO \cap \calV
  \label{defHOVO}
\Eeq
and endow them with their natural topologies as subspaces of $\calH$ and~$\calV$, respectively.
We observe that the first components $v$ of the elements $(v,\vG)\in\VO$
cannot span the whole of $C^\infty_c(\Omega)$
because of the zero mean value condition.
This has the following consequence:
variational equations with test functions in $\VO$
cannot be immediately read as equations in the sense of distributions
(this is \juerg{ the price we have to pay for the transformation of the old adjoint} system into the new~one!).
Hence, some care is \juerg{in order, and we have to prove some auxiliary lemmas. Here, we}
use the notation $\uG$ even though
it has nothing to do with the control variable.

\Blem
\label{Spanvz}
The set $\graffe{\vG:\ (v,\vG)\in\VO}$ is the whole of $\VG$.
\Elem

\Bdim
Take any $\vG\in\VG$.
As $\VG\subset\HxG{1/2}$,
there exists \juerg{some} \pier{$g\in\Huno$} such that $\pier{g} \suG=\vG$.
Now, we fix a closed ball $B\subset\Omega$ 
and \pier{a function  $\zeta \in C^1(\overline\Omega)$ such that
$\zeta =0$ in $\Omega\setminus B$  and $\int_B\zeta =1$.
Next, we define 
$m=\iO g $
and
$v:= g -m\zeta$.
Then, it turns out that $v \in\Huno$, 
$v\suG=g\suG = \vG$,}
and $\iO v=0$, i.e., $(v,\vG)\in\VO$.
\Edim

\Blem
\label{Interpr}
Assume that $(u,\uG)\in\calH$.
Then the condition
\Beq
  \iO uv + \iG \uG \vG = 0 
  \quad \hbox{for every $(v,\vG)\in\VO$}
  \label{hpinterpr}
\Eeq
implies that $u$ is a constant, \juerg{namely,} the mean value $\uO$ of~$u$, and $\uG=0$.
Moreover, $u=0$ if \eqref{hpinterpr} holds for every $(v,\vG)\in\calV$.
\Elem

\Bdim
We first decouple \eqref{hpinterpr}.
To this end, we fix $\vz\in\Hunoz$ such that $\vz^\Omega=1$
and set $k:=|\Omega|^{-1}\iO u\,\vz$.
Now, we take any $v\in\Hunoz$ and observe that
$v-\vO\vz$ belongs to $\Hunoz$ and has zero mean value.
Hence, $(v-\vO\vz,0)\in\VO$, and \eqref{hpinterpr} yields that
\Beq
  0 = \iO u (v-\vO\vz)
  = \iO u \, v - k \iO v
  = \iO (u - k) v.
  \non
\Eeq
As $v\in\Hunoz$ is arbitrary and $\Hunoz$ is dense in $\Ldue$,
we infer that $u=k$ \aeO, i.e., $u$ is a constant,
and  \juerg{this constant must equal} ~$\uO$.
Hence, \eqref{hpinterpr} implies
\Beq
  \iG \uG \vG = 0 
  \quad \hbox{for every $(v,\vG)\in\VO$}.
  \non
\Eeq
By Lemma~\ref{Spanvz}, the above equality holds for every $\vG\in\VG$.
As this space is dense in~$\LxG2$, we deduce that $\uG=0$.
If in addition \eqref{hpinterpr} holds for every $(v,\vG)\in\calV$,
then we can take $v=1$ and $\vG=1$ in \eqref{hpinterpr}
and deduce that $\uO=0$
(since we already know that $\uG=0$).
\Edim

\Bcor
\label{Density}
The space $\VO$ is dense in $\HO$.
\Ecor

\Bdim
We prove the following equivalent statement:
the only element $(u,\uG)\in\HO$ 
that is orthogonal to $\VO$
with respect to the scalar product in $\calH$ 
is the zero element of~$\HO$.
Thus, we assume that
\Beq
  \iO u \, v + \iG \uG \, \vG = 0 
  \quad \hbox{for every $(v,\vG)\in\VO$}.
  \label{hpdensity}
\Eeq
By Lemma~\ref{Interpr}, we deduce that
$u$ is a constant and that $\uG=0$.
As $u\in\HO$, the constant must be~$0$.
Therefore, $(u,\uG)=(0,0)$.
\Edim

In order to simplify the form of the problems we are dealing with, we introduce a notation.
Starting from the state $(\yopt,\yGopt)$
associated to an optimal control,
we~set
\Bsist
  && \lambda := f''(\yopt) \,, \quad
  \lamG := \fG''(\yGopt)
  \label{deflamlamG}
  \\
  && \phQ := \bQ (\yopt-\zQ) \,, \quad
  \phS := \bS (\yGopt-\zS) 
  \label{defphiQS}
  \\
  && \phO := \bO \bigl( \yopt(T) - \zO \bigr) \,, \quad
  \phG := \bG \bigl( \yGopt(T) - \zG \bigr) \,.
  \label{defphiOG}
\Esist
\gianni{Then, the adjoint problem \Pbladj\ becomes:}
\Bsist
  && \iO q(t) \, v
  = \iO \nabla p(t) \cdot \nabla v
  \quad \hbox{for \pier{all $t\in [0,T]$ and $v\in V$}}
  \label{primaV}
  \\
  && \gianni{- \iO \dt \bigl( p(t) + q(t) \bigr) \, v
  + \iO \nabla q(t) \cdot \nabla v
  + \iO \lambda(t) \, q(t) \, v}
  \qquad
  \non
  \\
  && \qquad {}
  - \iG \dt\qG(t) \, \vG
  + \iG \nablaG\qG(t) \cdot \nablaG\vG
  + \iG \lamG(t) \, \qG(t) \, \vG
  \non
  \\
  &&  \quad {}
  = \iO \phQ(t) v
  + \iG \phS(t) \vG
  \quad \hbox{\aat\ and every $(v,\vG)\in\gianni\calV$}
  \qquad
  \label{secondaV}
  \\
  && \iO (p+q)(T) \, v 
  + \iG \qG(T) \, \vG
  = \iO \phO v
  + \iG \phG \vG
  \quad \hbox{for every $(v,\vG)\in\gianni\calV$} .
  \qquad
  \label{cauchyV}
\Esist
\Accorpa\AggiuntoV primaV cauchyV
The result stated below ensures \juerg{the equivalence of problem \AggiuntoV\ 
and a new problem with decoupled equations}, as sketched at the beginning of the present section.
We note at once that the latter is plainly meaningful 
since $\calN q$ is well defined (see~\eqref{predefN}).
The statement also involves the operator $\calM:\L2{\Hx2}\to H^1(0,T)$  
defined~by
\Beq
  (\calM(v))(t) :=
  \phOO
  - \frac 1{|\Omega|} \int_t^T \!\! \iO \bigl( - \Delta v + \lambda v - \phQ \bigr) 
  \quad \hbox{for every $t\in[0,T]$}.
  \label{defM}
\Eeq
We notice that the subsequent proof will also show that the adjoint problem
is solved in the strong form presented in the Introduction.

\Bthm
\label{Equivalenza}
Assume \Regadj.
\gianni{Then, $(p,q,\qG)$ solves problem \AggiuntoV\ if and only~if}
\Beq
  q^\Omega(t) = 0
  \aand
  p(t) = \calN(q(t)) + (\calM(q))(t)
  \quad \pier{\hbox{for every } t\in [0,T]} 
  \label{primaN}
\Eeq
\begin{align}
  & - \iO \dt \bigl( \calN(q(t)) + q(t) \bigr) \, v
  + \iO \nabla q(t) \cdot \nabla v
  + \iO \lambda(t) \, q(t) \, v
  \qquad
  \non
  \\
  & \qquad {}
  - \iG \dt\qG(t) \, \vG
  + \iG \nablaG\qG(t) \cdot \nablaG\vG
  + \iG \lamG(t) \, \qG(t) \, \vG
  \non
  \\
  & \quad {}
  = \iO \phQ(t) v
  + \iG \phS(t) \vG
  \quad \hbox{\aat\ and every $(v,\vG)\in\VO$}
  \label{secondaN}
\end{align}
\Beq
  \iO \bigl( \calN q + q \bigr)(T) \, v 
  + \iG \qG(T) \, \vG
  = \iO \phO v
  + \iG \phG \vG
  \quad \hbox{for every $(v,\vG)\in\VO$} \,.
  \label{cauchyN}
\Eeq
\Accorpa\AggiuntoN primaN cauchyN
\Accorpa\Nuovaq secondaN cauchyN
\Ethm

\Bdim
We assume that \gianni{$(p,q,\qG)$ satisfies} \AggiuntoV\
and prove that \gianni{it} solves \AggiuntoN.
We often omit writing the time $t$ in order to simplify the notation.
By taking $v=1$ in~\eqref{primaV}, 
we see that the first \juerg{assertion of \eqref{primaN} holds.
In particular, the second assertion of (\ref{primaN}) is meaningful.}
\gianni{Moreover, by the definition \eqref{defN} of~$\calN$, we~have
\Beq
  p(t) - \pO(t) = \calN(q(t))
  \quad \hbox{or} \quad
  p(t) = \calN(q(t)) + \pO(t)
  \quad \aat .
  \label{reconstr}
\Eeq
}%
We now prove the second \pier{equality in}~\eqref{primaN}.
By taking any $v\in \pier{\calD (\Omega)}$ 
and using $(v,0)$ as a test functions in~\eqref{secondaV},
we derive~that
\Beq
  - \dt (p+q) - \Delta q + \lambda q = \phQ
  \quad \hbox{or} \quad
  \frac 1 {|\Omega|} \, \dt (p+q)
  = \frac 1 {|\Omega|} \, ( -\Delta q + \lambda q - \phQ )
  \non
\Eeq
in the sense of distributions on~$Q$, whence \aeQ\ as well,
due to the regularity of $p$ and~$q$.
By observing that both $q$ and $\dt q$ 
have zero mean value
(the~latter as a consequence of the former),
and just integrating the last equation over~$\Omega$,
we obtain 
\Beq
  \frac {d\pO}{dt} 
  = \frac 1 {|\Omega|}
  \iO ( -\Delta q + \lambda q - \phQ )
  \quad \hbox{whence} \quad
  \pO(t) = \pO(T) - \frac 1 {|\Omega|} \int_t^T \!\! \iO ( -\Delta q + \lambda q - \phQ ) \,.
  \non
\Eeq
On the other hand, \eqref{cauchyV} implies that $(p+q)(T)=\phO$,
whence $\pO(T)=\phOO$ since $q(T)$ has zero mean value.
By combining, we infer~that
\Beq
  \pO(t)
  = \phOO
  - \frac 1 {|\Omega|} \int_t^T \!\! \iO ( -\Delta q + \lambda q - \phQ )
  = (\calM(q))(t) \,.
  \non
\Eeq
Therefore, the second \juerg{assertion in} \eqref{primaN} follows from~\eqref{reconstr}.
In order to prove~\accorpa{secondaN}{cauchyN},
it suffices to write \accorpa{secondaV}{cauchyV} with $(v,\vG)\in\VO$,
\gianni{by recalling \eqref{reconstr} once more
and observing that $\dt\pO$ and $\pO(T)$ are space independent}.

Conversely, we \juerg{now} assume that $(p,q,\qG)$ solves \AggiuntoN\
and prove that the equations \AggiuntoV\ are satisfied.
We start from~\eqref{primaN}.
As $\calM(q)$ is space independent,
by recalling the definition \eqref{defN} of the operator~$\calN$, 
we~have, \aat\ and every $v\in V$,
\Beq
  \iO \nabla p(t) \cdot \nabla v
  = \iO \nabla\calN q(t) \cdot \nabla v
  = \iO q(t) v \,.
  \non
\Eeq
This is exactly~\eqref{primaV}.
Now, we prove~\eqref{secondaV}.
We deduce the strong form of the problem
hidden in the variational equation \eqref{secondaN}
\pier{thanks} to the integration by parts formulas
\Bsist
  && \iO \nabla q(t) \cdot \nabla v
  = \iO (-\Delta q(t)) v
  + \iG \dn q(t) \, v\suG
  = \iO (-\Delta q(t)) v
  + \iG \dn q(t) \, \vG
  \non
  \\
  && \iG \nablaG\qG(t) \cdot \nablaG\vG
  = \iG (-\DeltaG \pier{\qG} (t)) \, \vG ,
  \non
\Esist
where $(v,\vG)\in\calV$.
Thus, \eqref{secondaN} becomes
\Beq
  \iO u(t) v + \iG \uG(t) \vG = 0
  \quad \hbox{\aat\ and every $(v,\vG)\in\VO$,}
  \non
\Eeq
where the pair $(u,\uG)\in\L2\calH$ is given~by
\Bsist
  && u(t) := -\dt \bigl( \calN(q(t)) + q(t) \bigr)
  - \Delta q(t) + \lambda(t) \, q(t) - \phQ(t)
  \non
  \\
  && \uG(t) := \dn q(t) - \dt\qG(t) - \DeltaG\qG(t) + \lamG(t) \, \qG(t) - \phS(t) \,.
  \non
\Esist
Then, Lemma~\ref{Interpr} yields
\Beq
  u(t) = (u(t))^\Omega
  \aand
  \uG(t) = 0
  \quad \aat .
  \non
\Eeq
On the other hand, by recalling that $\calN(\dt q)$ and $\dt q$ have zero mean values
by the definition of $\calN$ and the first identity in \eqref{primaN},
and owing to the definition \eqref{defM} of~$\calM$, we~have
\Bsist
  && |\Omega| \, (u(t))^\Omega
  = \iO \bigl\{
  - \bigl( \calN(\dt q(t)) + \dt q(t) \bigr)
  - \Delta q(t) + \lambda(t) \, q(t) - \phQ(t)
  \bigr\}
  \non
  \\
  && = \iO \bigl\{
  - \Delta q(t) + \lambda(t) \, q(t) - \phQ(t)
  \bigr\}
  = |\Omega| \, \dt (\calM(q))(t)
  \quad \aat .
  \non
\Esist
We infer~that
\Beq
  -\dt \bigl( \calN(q(t)) + q(t) \bigr)
  - \Delta q(t) + \lambda(t) \, q(t) - \phQ(t)
  = \dt (\calM(q))(t)
  \quad \aat
  \non
\Eeq
so that \eqref{primaN} yields
\Bsist
  && \iO \bigl\{
  - \dt \bigl( p(t) + q(t) \bigr)
  - \Delta q(t) + \lambda(t) \, q(t) - \phQ(t)
  \bigr\} \, v
  = 0
  \non
  \\
  && \quad \hbox{\aat\ and every $v\in V$}.
  \non
\Esist
Now, the identity $\uG=0$ implies
\Bsist
  && \iG \bigl\{
  \dn q(t) - \dt\qG(t) - \DeltaG\qG(t) + \lamG(t) \, \qG(t) - \phS(t)
  \bigr\}\, \pier{\vG}
  = 0
  \non
  \\
  && \quad \hbox{\aat\ and every $\vG\in\VG$}.
  \non
\Esist
In particular, for any $(v,\vG)\in\calV$, we can write
both previous equalities and add them to each other.
By integrating by parts in the opposite direction,
we deduce~\eqref{secondaV}.
Finally, by applying Lemma~\ref{Interpr} once more,
we derive from~\eqref{cauchyN}
\Beq
  \bigl( \calN q + q \bigr)(T) - \phO = k
  \aand
  \qG(T) - \phG = 0
  \non
\Eeq
where $k$ is the mean value of the \lhs.
\pier{Note that} both $q(T)$ and $(\calN q)(T)$ have zero mean value
by the first identity in \eqref{primaN} and the definition \eqref{defN} of~$\calN$\pier{. Then},
we have $k=-\phOO$.
Hence, by the definition of~$\calM$, we infer~that
\Beq
  (\calN q+q)(T)
  = \phO - \phOO
  = \phO - (M(q))(T) .
  \non
\Eeq
Then, the second assertion in \eqref{primaN} yields $(p+q)(T)=\phO$,
and \eqref{cauchyV} follows immediately.\Edim

Thanks to the theorem just proved, we can replace the old adjoint problem
by the new one \juerg{in which the equations are decoupled}.
As we are going to see the sub-problem for $(q,\qG)$ 
as an abstract differential equation,
we prepare the proper framework, which is related to the Hilbert spaces
$\VO$ and $\HO$ defined in~\eqref{defHOVO}. \juerg{To this end, let
$_{{\cal V}^*_\Omega} \langle \,\cdot\,,\,\cdot\,\rangle_{{\cal V}_\Omega}\,$ denote  the
dual pairing between ${\cal V}^*_\Omega$ and ${\cal V}_\Omega$.
Then, recalling that $\VO$ is by Corollary~\ref{Density} dense in~$\HO$,
we can construct the Hilbert triplet $(\VO,\HO,\VOp)$, that is,
we identify $\HO$ with a subspace of~$\VOp$, the dual space of~$\VO$,
in order that
\begin{equation}
\label{ternahilb}
_{\VOp} \< (u,\uG) , (v,\vG) > _{\VO}
  = \bigl( (u,\uG) , (v,\vG) \bigr)_{\HO}
\quad\forall\,(u,\uG)\in\HO, \quad\forall\,(v,\vG)\in\VO\,.
\end{equation}
Here, we define the scalar product $(\cpto,\cpto)_{\HO}$ and the scalar product in~$\VO$ 
by}
\Bsist
  & \bigl( (u,\uG) , (v,\vG) \bigr)_{\HO}
  & := \iO u \, v + \iG \uG \, \vG 
  \label{scalHO}
  \\
  & \bigl( (u,\uG) , (v,\vG) \bigr)_{\VO}
  & := \iO \nabla u \cdot \nabla v + \iG \nablaG\uG \cdot \nablaG\vG \,.
  \label{scalVO}
\Esist
In \eqref{scalHO} (resp.~\eqref{scalVO}), 
$(u,\uG)$ and $(v,\vG)$ \pier{denote} generic elements of~$\HO$ (resp.~$\VO$).
Note that \eqref{scalVO} actually defines a scalar product in $\VO$
that is equivalent to the standard one by the Poincar\'e inequality~\eqref{poincare}.
We also introduce the associated Riesz operator $\RO\in\calL(\VO,\VOp)$, namely
\Beq
  _{\VOp} \< \RO (u,\uG) , (v,\vG) > _{\VO}
  = \bigl( (u,\uG) , (v,\vG) \bigr)_{\VO}
  \quad \hbox{for every $(u,\uG),(v,\vG)\in\VO$} \,.
  \label{riesz}
\Eeq

\juerg{Since, as already mentioned, variational equations with test functions in $\VO$
cannot immediately be read as differential equations,}
we also prove the following lemma.

\Blem
\label{Perregolarita}
Assume $(u,\uG)\in\VO$ and $\RO(u,\uG)\in\HO$.
Then we have $u\in\Hx2$ and $\uG\in\HxG2$.
Moreover, it holds
\Beq
  \norma u_{\Hx2} + \norma\uG_{\HxG2}
  \leq c\, \norma{\RO(u,\uG)}_{\HO}\,, 
  \label{perregolarita}
\Eeq
where $c$ depends only on~$\Omega$.
\Elem

\Bdim
The assumptions mean that there exists some $(\psi,\psiG)\in\HO$ such that
\Bsist
  && \bigl( (u,\uG) , (v,\vG) \bigr)_{\VO}
  = \bigl( (\psi,\psiG) , (v,\vG) \bigr)_{\HO} \,,
  \quad \hbox{that is,}
  \non
  \\
  && \iO \nabla u \cdot \nabla v + \iG \nablaG\uG \cdot \nablaG\vG
  = \iO \psi v + \iG \psiG \vG
  \label{hpreg}
\Esist
for every $(v,\vG)\in\VO$.
As in the proof of Lemma~\ref{Interpr}, we decouple~\eqref{hpreg}.
We fix $\vz\in\Hunoz$ such that $\vz^\Omega=1$ and set
\juerg{
$$c_1:=\iO\nabla u\cdot\nabla\vz, \quad
c_2:=\iO\psi\vz$$
}
and $k_i:=c_i/|\Omega|$ for $i=1,2$.
Now, we take any $v\in\Hunoz$.
As $v-\vO\vz$ \juerg{belongs} to $\Hunoz$ and has zero mean value,
we have $(v-\vO\vz,0)\in\VO$, and \eqref{hpreg} yields
\Beq
  \iO \nabla u \cdot \nabla(v-\vO\vz)
  = \iO \psi (v-\vO\vz)
  \quad \hbox{or} \quad
  \iO (\nabla u \cdot \nabla v - k_1 v)
  = \iO (\psi - k_2) v \,.
  \non
\Eeq
As $v\in\Hunoz$ is arbitray, this simply means
\Beq
  - \Delta u = \psi + k
  \quad \hbox{where $k:=k_1-k_2$}.
  \label{perregO}
\Eeq
In particular, \pier{we infer that $\Delta u\in\Ldue$ and this, combined with $u\suG = \uG \in H^1(\Gamma)$, yields (cf., e.g., \cite[Thm.~3.2, p.~1.79]{BreGil}) $u\in H^{3/2} (\Omega)$. Then, 
by a trace theorem stated, e.g., in \cite[Thm.~2.25, p.~1.62]{BreGil} it follows that 
$\dn u$ lies in $\LxG2$ and we can integrate by parts.
Hence, for any $(v,\vG)\in\VO$ we have}
\Bsist
  && \iO (\psi + k) v
  + \iG \nablaG\uG \cdot \nablaG\vG
  = \iO (-\Delta u) v
  + \iG \nablaG\uG \cdot \nablaG\vG
  \non
  \\
  && = \iO \nabla u \cdot \nabla v
  \pier{{}- \iG \dn u \, v\suG} 
  + \iG \nablaG\uG \cdot \nablaG\vG
  \non
  \\
  && = \iO \psi v
  + \iG \psiG \vG 
  \pier{{}- \iG \dn u \, \vG}
  = \iO (\psi + k) v
  + \iG (\psiG  
  \pier{{}-  \dn u ) \vG .}
  \non
\Esist
\pier{Therefore, we deduce that}
\Beq
  \iG \nablaG\uG \cdot \nablaG\vG
  \pier{{}= \iG (\psiG - \dn u ) \vG}
   \quad \hbox{for every $(v,\vG)\in\VO$}
  \non
\Eeq
and \pier{Lemma~\ref{Spanvz} implies that 
the same equality holds for every $\vG\in\VG$,
whence} 
\Beq
  - \DeltaG\uG  = \psiG - \dn u  \quad \hbox{on $\Gamma$}.
  \label{regDeltaG}
\Eeq
\pier{As $\psiG$ and $\dn u$ both belong to $\LxG2$,
the regularity theory for elliptic equation
(in~fact, its boundary version)
implies that $\uG\in\HxG{2}$ (see, e.g., \cite[Thms.~7.4 and~7.3, pp.\ 187-188]{LioMag}
or \cite[Thm.~3.2, p.~1.79, and Thm.~2.27, p.~1.64]{BreGil}).}
Coming back to~$u$, we thus have $\Delta u\in\Lx2$ and $u\suG=\uG\in\HxG2$,
whence $u\in\Hx2$.
Finally, as each of the regularity results we have used
corresponds to an estimate for the related norm,
\eqref{perregolarita}~holds as well.
\Edim

\vspace{5mm}
At this point, we are ready to prove Theorem~\ref{Existenceadj}.
Thanks to Theorem~\ref{Equivalenza},
it is sufficient to prove that there exists a unique solution to problem~\AggiuntoN\
satisfying the regularity requirements~\Regadj.
Moreover, once the existence of a unique solution $(q,\qG)$ to \Nuovaq\
with the prescribed regularity is established,
it suffices to observe that \eqref{primaN} provides
a function $p$ that fulfills~\eqref{regp}.
Indeed, \eqref{regq} and \eqref{regN} imply $\calN q\in\L2{\Hx4}$
and $\dt\calN q=\calN(\dt q)\in\L2{\Hx2}$.
On the other hand, $\calM(q)$ is space independent, and its time derivative belongs to $L^2(0,T)$
since it is the mean value of an element of~$\LQ2$.
Hence,  \pier{we have that $p\in\L2{\Hx4}$ and $\dt p\in\L2{\Hx2}$.}

\juerg{In the following, we denote pairs belonging to ${\cal H}_\Omega$ by bold letters, writing, for instance,
${\bf v}$ in place of $(v,v_\Gamma)$. From this no confusion will arise.} 
We are going to present the problem in the form
\Bsist
  && - \frac d{dt} \, \bigl( \calB\qqq(t),\vvv \bigr)_{\HO} 
  + {}{}_{\VOp} \< \calA(t)\qqq(t) , \vvv >_{\VO}
  = {}{}_{\VOp} \< \fff(t),\vvv >_{\VO}
  \non
  \\
  && \quad \hbox{\aat\ and every $\vvv\in\VO$}
  \label{preastratta}
  \\
  && \bigl( (\calB\qqq)(T) , \vvv \bigr)_{\HO}
  = \bigl( \zT , \vvv \bigr)_{\HO}
  \quad \hbox{for every $\vvv\in\HO$} 
  \label{precauchyastratto}
\Esist
with a proper choice of the operators
$\calA(t)\in\calL(\VO,\VOp)$ and $\calB\in\calL(\VO,\HO)$, 
and of the data $\fff\in\L2\VOp$ and $\zT\in\HO$.
This means \pier{that} the following backward Cauchy problem
\Beq
  -\frac d{dt} \, \bigl( \calB \, \qqq(t) \bigr) + \calA(t) \, \qqq(t) = \fff(t)
  \quad \aat,
  \aand
  (\calB\qqq)(T) = \zT 
  \label{pblastratto}
\Eeq
\Accorpa\Preastratto preastratta precauchyastratto
\pier{has to be solved.} Problem \eqref{pblastratto} 
(in~fact the equivalent forward problem obtained 
by replacing $t$ by~$T-t$)
is~well known
(see \cite{Baio} for a very general situation that allows for time dependent and even nonlocal operators).
Here, we recall sufficient conditions that imply those given in  
\cite[Thm.~7.1, p.~70]{Lions}
and thus yield well-posedness in a proper framework. 
We can require~that
\Bsist
  && \calA(t) = \Az + \Lambda(t) \,, \quad
  \Az \in \calL(\VO,\VOp)
  \aand
  \Lambda(t) \in \calL(\HO,\HO)
  \quad \aat;
  \non
  \\[2mm]
  && _{\VOp} \< \Az \vvv , \vvv > _{\VO}
  \geq \alpha \norma\vvv_{\VO}^2
  \quad \hbox{for some $\alpha>0$ and every $\vvv\in\VO$}; 
  \non
  \\[2mm]
  && \norma{\Lambda(t)\vvv}_{\HO}
  \leq M \norma\vvv_{\HO}
  \quad \hbox{for some constant $M>0$ and every $\vvv\in\HO$}; 
  \non
  \\[2mm]
  && \calB \in \calL(\HO,\HO)
  \quad \hbox{is symmetric and satisfies}
  \non
  \\
  && _{\HO} \bigl( \calB \vvv , \vvv \bigr) _{\HO}
  \geq \alpha \norma\vvv_{\HO}^2
  \quad \hbox{for some $\alpha>0$ and every $\vvv\in\VO$} \,. 
  \non
\Esist
Moreover, $\Lambda$ is (properly) measurable with respect to~$t$.
If such conditions hold then,
for every $\fff\in\L2\VOp$ and $\zT\in\HO$,
problem \eqref{pblastratto} has a unique solution
\Beq
  \qqq \in \H1\VOp \cap \L2\VO \subset \C0\HO \,.
  \non
\Eeq
Furthermore, the solution $\qqq$ also satisfies
\Beq
  \qqq' \in \L2\HO
  \quad \hbox{if} \quad
  \hbox{$\Az$ is symmetric, \ $\fff\in\L2\HO$ \ and \ $\qqq(T)\in\VO$}.
  \non
\Eeq
In our case, \pier{we choose}
\Bsist
  && \Az = \RO \,, 
  \quad \hbox{the Riesz operator \eqref{riesz}}
  \non
  \\[1mm]
  && \Lambda(t) (v,\vG) = \bigl( \lambda(t) v - (\lambda(t) v)^\Omega , \lamG(t) \vG \bigr) 
  \quad \hbox{\aat\ and $(v,\vG)\in\HO$}
  \non
  \\[1mm] 
  && \calB (v,\vG) = (\calN v + v , \vG) 
  \quad \hbox{for every $(v,\vG)\in\HO$}
  \non
  \\[1mm]
  && \fff(t) := \bigl( \phQ(t) - (\phQ(t))^\Omega , \phS(t) \bigr) 
  \quad \aat
  \non
  \\[1mm]
  && \zT := \bigl( \phO - \phOO , \phG \bigr) \,.
  \non
\Esist
The choices of $\Az$ and $\calB$ being clear, 
the other ones exactly yield what we need,~i.e.,
\Bsist
  && \bigl( \Lambda(t) (u,\uG) , (v,\vG) \bigr)_{\HO}
  = \iO \lambda(t) \, uv
  + \iG \lamG(t) \, \uG \vG
  \non
  \\
  && \quad \hbox{\aat\ and $(u,\uG),\,(v,\vG)\in\HO$}
  \non
  \\
  && _{\VOp} \< \fff(t),(v,\vG) >_{\VO}
  = \iO \phQ(t) v + \iG \phS(t) \vG 
  \non
  \\
  && \quad \hbox{\aat\ and $(v,\vG)\in\VO$} 
  \non
  \\
  && \bigl( \zT , (v,\vG) \bigr)_{\HO}
  = \iO \phO v + \iG \phG \vG 
  \quad \hbox{for $(v,\vG)\in\HO$} \,.
  \non
\Esist
Furthermore, the conditions we have required on the operators are fulfilled.
Indeed,
\Beq
  \normaH{\lambda(t) v - (\lambda(t) v)^\Omega}
  + \normaHG{\lamG(t) \vG}
  \leq c \bigl( \normaH v + \normaHG\vG \bigr)
  \non
\Eeq
\aat\ and every $(v,\vG)\in\HO$,
since the functions $\lambda$ and $\lamG$ are bounded (see~\eqref{deflamlamG}),
and $\calB$ is \juerg{symmetric} and coercive
since $\calN$ is symmetric and positive (see, in particular,~\eqref{simmN}).
Finally, by accounting for \eqref{hpzzzz}, \accorpa{defphiQS}{defphiOG}, \eqref{hpzz},
we see that $\fff\in\L2\HO$ and $\qqq(T)=(0,0)$.
Therefore, problem \Nuovaq\ has a unique solution satisfying
\Beq
  (q,\qG) \in \H1{\HO} \cap \L2\VO 
  \aand
  \RO(q,\qG) \in \L2{\HO},
  \non
\Eeq
the last one by comparison in~\eqref{pblastratto}.
Then, Lemma~\ref{Perregolarita} ensures that
$q(t)\in\Hx2$ and $\qG(t)\in\HxG2$ \aat\ and that the estimate
\Beq
  \norma{q(t)}_{\Hx2} + \norma{\qG(t)}_{\HxG2}
  \leq c \norma{\RO(q(t),\qG(t))}_{\HO} 
  \non
\Eeq
holds true \aat.
This implies $q\in\L2{\Hx2}$ and $\qG\in\L2{\HxG2}$, and the proof is complete.
\juerg{\qed}

\Brem
\label{Pesi}
Assumption \eqref{hpzz} can be avoided provided that
we require less regularity from the solution to the adjoint problem.
\pier{More precisely, we keep the regularity related to the variational structure of 
\accorpa{primaadj}{cauchyadj}, e.g., we ask for $q\in\H1{\Vp} \cap \L2V$, while we replace $L^2$ 
summability by weighted $L^2$ summability where spaces of smooth functions on $\Omega$ are involved, e.g., $q\in \L2{H^2(\Omega)}$. Namely, we require~that}
\Bsist
  && \hbox{the functions} \quad
  t\mapsto (T-t)^{1/2}q(t) 
  \aand
  t\mapsto (T-t)^{1/2}\dt q(t) 
  \non
  \\
  &&  \hbox{belong to $\L2{\Hx2}$ and to $\L2H$, respectively}
  \non
\Esist
and we analogously \juerg{deal} with the other conditions.
By doing that, the equivalence stated in Theorem~\ref{Equivalenza} still holds.
On the other hand, the derivative $\qqq'$ of the solution $\qqq$ to the abstract problem
satisfies the right weighted summability that yields the new requirements
without assuming that $\qqq(T)\in\VO$,
so that \eqref{hpzz} is not needed.
For the reader's convenience, we sketch the formal a~priori estimate
that yields the mentioned property of~$\qqq'$
whenever it is replaced by a rigorous argument.
For convenience, we set 
\Beq
  u(t):=\qqq(T-t), \quad
  \mu(t) := \Lambda(T-t)
  \aand
  g(t) := \fff(T-t)
  \non
\Eeq
and write \eqref{pblastratto} as a forward Cauchy problem for~$u$.
Then, we formally test the new equation by $tu'(t)$
and integrate with respect to time.
We simply write $(\cpto,\cpto)$ for both the duality pairing between $\VOp$ and $\VO$
and for the scalar product in~$\HO$.
We have, for every $t\in[0,T]$,
\Beq
  \iot \bigl( \calB u'(s) , s u'(s) \bigr) \, ds
  + \iot \bigl(\Az u(s) , s u'(s) \bigr) \, ds
  = \iot \bigl( g(s) - \mu(s) u(s), s u'(s) \bigr) \, ds \,.
  \non
\Eeq
As $\Az$ is symmetric and both $\Az$ and $\calB$ are coercive, 
we can estimate the \lhs\ from below as follows
\Bsist
  && \iot \bigl( \calB u'(s) , s u'(s) \bigr) \, ds
  + \iot \bigl( \Az u(s) , s u'(s) \bigr) \, ds
  \non
  \\
  && 
  = \iot s \bigl( \calB u'(s) , u'(s) \bigr) \, ds
  + \frac 12 \iot \frac d{ds} \, \bigl\{ s \bigr( \Az u(s), u(s) \bigr) \bigr\} \, ds
  - \frac 12 \iot \bigl( \Az u(s) , u(s) \bigr) \, ds
  \non
  \\
  && \geq \alpha \iot s \norma{u'(s)}_{\HO}^2 \, ds
  + \frac \alpha 2 \, t \norma{u(t)}_{\VO}^2 
  - c \ioT \norma{u(s)}_{\VO}^2 \, ds \,.
  \non
\Esist
On the other hand, as $\norma{\mu(t)}_{\calL(\HO,\HO)}\leq M$ \aat, we also have
\Bsist
  && \iot \bigl( g(s) - \mu(s) u(s) , s u'(s) \bigr) \, ds
  \non
  \\
  && \leq \frac \alpha 2 \iot s \norma{u'(s)}_{\HO}^2 \, ds
  + c \ioT s \norma{g(s)}_{\HO}^2 \, ds
  + c \ioT s \norma{u(s)}_{\HO}^2 \, ds \,.
  \non
\Esist
By combining we infer that
\Beq
  \iot s \norma{u'(s)}_{\HO}^2 \, ds
  + t \norma{u(t)}_{\VO}^2 
  \leq c \bigl( \norma g_{\L2\HO}^2 + \norma u_{\L2\VO}^2 \bigr) 
  \quad \hbox{for every $t\in[0,T]$}.
  \non
\Eeq
As the last norm of $u$ is supposed to be already estimated, 
we obtain the desired weighted summability for~$u'$
as well as a weighted boundedness for $u$ in~$\VO$, as a by-product.
\Erem


\section{Necessary optimality conditions}
\label{OPTIMALITY}
\setcounter{equation}{0}

In this section, we derive the optimality condition~\eqref{cnoptadj}
stated in Theorem~\ref{CNoptadj}.
We start from~\eqref{precnopt} and first prove~\eqref{cnopt}.
We recall the definitions \Controltostate\ of the spaces $\calX$ and $\calY$
and of the control-to-state mapping~$\calS$.

\Bprop
\label{CNopt}
Let $\uopt$ be an optimal control and $(\yopt,\yGopt):=\calS(\uopt)$.
Then, \eqref{cnopt} holds.
\Eprop

\Bdim
This is just due to the chain rule for \Frechet\ derivatives, 
as already said in Section~\ref{STATEMENT},
and we just provide some detail.
Let $\tilde\calS:\calU\to\calY\times\calX$ be given by
$\tilde\calS(\uG):=(\calS(\uG),\uG)$.
Then, $\tilde\calS$~is \Frechet\ differentiable at any $\uG\in\calU$ since $\calS$ is~so.
Precisely, thanks to Theorem~\ref{Fdiff},
the \Frechet\ derivative $D\tilde\calS(\uG)$ acts as follows
\Beq
  D\tilde\calS(\uG):
  \hG \mapsto \bigl( [D\calS(\uG)](\hG),\hG \bigr)
  = (\xi,\xiG,\hG)
  \quad \hbox{for $\hG\in\calX$}
  \non
\Eeq
where $(\xi,\xiG,\eta)$ is the solution to the linearized problem \Linpbl\ corresponding to~$\hG$.
On the other hand, if we see the cost functional \eqref{defcost}
as a map from $\calY\times\calX$ to~$\erre$,
it is clear that its \Frechet\ derivative $D\calJ(y,\yG,\uG)$ at $(y,\yG,\uG)\in\calY\times\calX$
is given~by
\Bsist
  && [D\calJ(y,\yG,\uG)](k,k_\Gamma,\hG)
  = \bQ \intQ (y-\zQ) k 
  + \bS \intS (\yG-\zS) k_\Gamma
  \non
  \\
  && \quad {}
  + \bO \iO (y(T)-\zO) k(T)
  + \bG \iG (y(T)-\zG) k_\Gamma(T)
  + \bz \intS \uG \hG 
  \non
  \\
  && \hbox{for $(k,k_\Gamma)\in\calY$ and $\hG\in\calX$}.
  \non
\Esist
Therefore, being $\redJ=\calJ\circ\tilde\calS$, 
the chain rule imples that $[D\redJ(\uG)]$ maps any $\hG\in\calX$ into
\Bsist
  && [D\redJ(\uG)](\hG)
  = [D\calJ(\tilde\calS(\uG)] \bigl( [D\tilde\calS(\uG)] (\hG) \bigr)
  \non
  \\
  && = [D\calJ(\tilde\calS(\uG)] (\xi,\xiG,\hG)
  = [D\calJ(y,\yG,\uG)] (\xi,\xiG,\hG)
  \non
  \\
  && = \bQ \intQ (y-\zQ) \xi 
  + \bS \intS (\yG-\zS) \xiG
  \non
  \\
  && \quad {}
  + \bO \iO (y(T)-\zO) \xi(T)
  + \bG \iG (y(T)-\zG) \xiG(T)
  + \bz \intS \uG \hG 
  \non
  \non
\Esist
where $(y,\yG)=\calS(\uG)$ and $(\xi,\xiG)$ has the same meaning as before.
Therefore, \eqref{cnopt} immediately follows from~\eqref{precnopt}.
\Edim

At this point, we are ready to prove Theorem~\ref{CNoptadj} on optimality,
i.e., the necessary condition \eqref{cnoptadj} for
$\uopt$ to be an optimal control in terms of the solution $(p,q,\qG)$
of the adjoint problem~\Pbladj.
We note that it is sufficient to prove the following:
if $\uG\in\calU$, $(y,\yG)=\calS(\uG)$, $\hG\in\calX$,
$(\xi,\xiG,\eta)$ is the solution to the linearized problem \Linpbl\ corresponding to~$\hG$,
and $(p,q,\qG)$ solves the adjoint problem \Pbladj\ (where one reads
$(y,\yG)$ in place of $(\yopt,\yGopt)$),
then
\Bsist
  & \displaystyle\intS \qG \hG 
  & = \intQ \bQ (y - \zQ) \xi
  + \intS \bS (\yG - \zS) \xiG
  \non
  \\
  && {}
  + \iO \bO (y(T) - \zO) \xi(T)
  + \iG \bG (\yG(T) - \zG) \xiG(T) \,.
  \label{tesiCN}
\Esist
Indeed, once this is proved,
we can apply it to \juerg{any optimal control} $\uG:=\uopt$,
and \pier{\eqref{cnoptadj}} follows from the necessary condition \eqref{cnopt}
already established in Proposition~\ref{CNopt}.
So, we fix $\uG\in\calU$ and $\hG\in\calX$, 
and write both the linearized problem
and the adjoint problem we are interested in, for the reader's convenience.
All the variational equations hold \aat,
but we avoid writing the time~$t$, for brevity.
We~have
\Bsist
  && \iO \dt\xi \, v 
  + \iO \nabla\eta \cdot \nabla v = 0
  \label{primaL}
  \\
  && \iO q \, v
  = \iO \nabla p \cdot \nabla v
  \label{primaA}
  \\
  \noalign{\allowbreak}
  && \iO \eta v
  = \iO \dt\xi \, v
  + \iG \dt\xiG \, v
  + \iO \nabla\xi \cdot \nabla v
  + \iG \nablaG\xiG \cdot \nablaG\vG
  \qquad
  \non
  \\
  && \quad {}
  + \iO \lambda \, \xi \, v
  + \iG \bigl( \lamG \, \xiG - \hG \bigr) \, \vG
  \label{secondaL}
  \\
  \noalign{\allowbreak}
  && \gianni{%
  - \iO \dt \bigl( p + q \bigr) \, v
  + \iO \nabla q \cdot \nabla v
  + \iO \lambda \, q \, v }
  \non
  \\
  && \quad {}
  - \iG \dt\qG \, \vG
  + \iG \nablaG\qG \cdot \nablaG\vG
  + \iG \lamG \, \qG \, \vG
  \qquad
  \non
  \\
  && = \iO \bQ \bigl( y - \zQ \bigr) v
  + \iG \bS \bigl( \yG - \zS \bigr) \vG \,.
  \label{secondaA}
\Esist
In the above equations, $\lambda:=f''(y)$ and $\lamG:=\fG''(\yG)$.
Moreover, \accorpa{primaL}{primaA} hold for every $v\in V$,
while \accorpa{secondaL}{secondaA} are satisfied for every $(v,\vG)\in\calV$.
Furthermore, $\xi(0)=0$ and
\Beq
  \iO (p+q)(T) \, v 
  + \iG \qG(T) \, \vG
  = \iO \bO \bigl( y(T) - \zO \bigr) v
  + \iG \bG \bigl( \yG(T) - \zG \bigr) \vG
  \label{cauchyA}
\Eeq
for every $(v,\vG)\in\calV$.
We choose $v=p$ in~\eqref{primaL},
$v=\eta$ in~\eqref{primaA}, 
$(v,\vG)=(q,\qG)$ in~\eqref{secondaL} 
and $(v,\vG)=-(\xi,\xiG)$ in~\eqref{secondaA}.
\gianni{Then}, 
by integrating over $(0,T)$ all the equalities we obtain, we~have
\Bsist
  && \intQ \dt\xi \, p
  + \intQ \nabla\eta \cdot \nabla p = 0
  \non 
  \\
  && \intQ q \, \eta
  = \intQ \nabla p \cdot \nabla\eta
  \non 
  \\
  \noalign{\allowbreak}
  && \intQ \dt\xi \, q
  + \intS \dt\xiG \, \qG
  + \intQ \nabla\xi \cdot \nabla q
  + \intS \nablaG\xiG \cdot \nablaG\qG
  \qquad
  \non
  \\
  && \quad {}
  + \intQ \lambda \, \xi \, q
  + \intS \bigl( \lamG \, \xiG - \hG \bigr) \, \qG  
  = \intQ \eta q  
  \non 
  \\
  \noalign{\allowbreak}
  && \intQ \dt \bigl( p + q \bigr) \, \xi
  - \intQ \nabla q \cdot \nabla\xi
  - \intQ \lambda \, q \, \xi
  \non
  \\
  && \quad {}
  + \intS \dt\qG \, \xiG
  - \intS \nablaG\qG \cdot \nablaG\xiG
  - \intS \lamG \, \qG \, \xiG
  \qquad
  \non 
  \\
  && = - \intQ \bQ \bigl( y - \zQ \bigr) \xi
  - \intS \bS \bigl( \yG - \zS \bigr) \xiG \,.
  \non
\Esist
At this point, we add the above equalities to each other and just simplify.
We obtain
\Bsist
  && \intQ \dt\xi \, (p+q)
  + \intQ \dt \bigl( p + q \bigr) \, \xi
  + \intS \dt\xiG \, \qG
  + \intS \dt\qG \, \xiG
  - \intS \hG \, \qG
  \non
  \\
  && = - \intQ \bQ \bigl( y - \zQ \bigr) \xi
  - \intS \bS \bigl( \yG - \zS \bigr) \xiG \,.
  \non
\Esist
By accounting for the Cauchy condition $\xi(0)=0$,
we can write an equivalent form as follows
\Beq
  \iO (p+q)(T) \, \xi(T)
  + \iG \qG(T) \, \xiG(T)
  = \intS \hG \, \qG
  - \intQ \bQ \bigl( y - \zQ \bigr) \xi
  - \intS \bS \bigl( \yG - \zS \bigr) \xiG \,.
  \non
\Eeq
At this point, we choose $(v,\vG)=(\xi(T),\xiG(T))$ in \eqref{cauchyA}
and get
\Beq
  \iO (p+q)(T) \, \xi(T)
  + \iG \qG(T) \, \xiG(T)
  = \iO \bO \bigl( y(T) - \zO \bigr) \xi(T)
  + \iG \bG \bigl( \yG(T) - \zG \bigr) \xiG(T) \,.
  \non
\Eeq
By comparison, we conclude that \eqref{tesiCN} holds.
This completes the proof of Theorem~\ref{CNoptadj}.\juerg{\qed}



\vspace{3truemm}


\Begin{thebibliography}{10}

\bibitem{Baio} 
\juerg{C. Baiocchi,}
Sulle equazioni differenziali astratte lineari del primo 
e del secondo ordine negli spazi di Hilbert\pier{,}
{\it Ann. Mat. Pura Appl. (4)\/} {\bf 76} (1967) 233-304.



\bibitem{BreGil}
F. Brezzi and G. Gilardi,
Part~1:
Chapt.~2, Functional spaces, 
Chapt.~3, Partial differential equations,
in ``Finite element handbook'',
H. Kardestuncer and D. \pier{H.\ Norrie}  eds.,
McGraw-Hill Book Company,
NewYork,
1987.

\bibitem{CahH} 
J. W. Cahn and J. E. Hilliard, 
Free energy of a nonuniform system I. Interfacial free energy, 
{\it J. Chem. Phys.\/}
{\bf 2} (1958) 258-267.

\bibitem{CaCo}
L. Calatroni and P. Colli,
Global solution to the Allen--Cahn equation with singular potentials and dynamic boundary conditions,
{\it Nonlinear Anal.\/} 
{\bf 79} (2013) 12-27.

\bibitem{CFP} 
R. Chill, E. Fa\v sangov\'a and J. Pr\"uss,
Convergence to steady states of solutions of the Cahn--Hilliard equation with dynamic boundary conditions,
{\it Math. Nachr.\/} 
{\bf 279} (2006) 1448-1462.

\bibitem{CFS}
P. Colli, M. H. Farshbaf-Shaker and J. Sprekels,
A deep quench approach to the optimal control of an Allen--Cahn equation
with dynamic boundary condition\pier{s} and double \pier{obstacles}, 
\pier{{\it Appl. Math. Optim.\/}, doi:10.1007/s00245-014-9250-8
(see also WIAS Preprint No. 1838 (2013), pp. 1-23)}.


\bibitem{CF} 
\pier{P. Colli and T. Fukao,
The Allen-Cahn equation with dynamic boundary conditions and mass constraints,
\lastrev{{\it Math Methods Appl Sci\/}, doi:10.1002/mma.3329} 
(see also preprint arXiv:1405.0116~[math.AP] (2014), pp.~1-23).}

\bibitem{CGPS}
\juerg{P. Colli, G. Gilardi, P. Podio-Guidugli and J. Sprekels, 
Distributed optimal control of a nonstandard system of phase field equations, {\it Cont. Mech. Thermodyn.\/}
{\bf 24} (2012) 437-459.}  

\bibitem{CGScon}
\juerg{P. Colli, G. Gilardi and J. Sprekels,
Analysis and boundary control of a nonstandard system of phase field equations, 
{\it Milan J. Math.\/} {\bf 80} (2012) 119-149.}

\bibitem{CGS}
P. Colli, G. Gilardi and J. Sprekels,
On the Cahn--Hilliard equation with dynamic 
boundary conditions and a dominating boundary potential,
\pier{{\it J. Math. Anal. Appl.\/} {\bf 419} (2014) 972-994.}

\bibitem{CS}
P. Colli and J. Sprekels,
Optimal control of an Allen--Cahn equation 
with singular potentials and dynamic boundary condition,
\lastrev{{\it SIAM J Control Optim.\/} {\bf 53} (2015) 213-234.}

\bibitem{EllSt} 
C. M. Elliott and A. M. Stuart, 
Viscous Cahn--Hilliard equation. II. Analysis, 
{\it J. Differential Equations\/} 
{\bf 128} (1996) 387-414.

\bibitem{EllSh} 
C. M. Elliott and S. Zheng, 
On the Cahn--Hilliard equation, 
{\it Arch. Rational Mech. Anal.\/} 
{\bf 96} (1986) 339-357.


\bibitem{GiMiSchi} 
G. Gilardi, A. Miranville and G. Schimperna,
On the Cahn--Hilliard equation with irregular potentials and dynamic boundary conditions,
{\it Commun. Pure Appl. Anal.\/} 
{\bf 8} (2009) 881-912.



\bibitem{HW1}
\juerg{M. Hinterm\"uller and D. Wegner, Distributed optimal control of the 
Cahn--Hilliard representsystem
including the case of a double-obstacle homogeneous free energy density, {\it SIAM J.
Control Optim.} {\bf 50} (2012) 388-418.}

\bibitem{HW2}
\juerg{M. Hinterm\"uller and D. Wegner,	Optimal control of a semi-discrete 
Cahn--Hilliard--Navier--Stokes system,} \pier{{\it SIAM J. Control Optim.} 
{\bf 52} (2014) 747-772.}

\bibitem{HK}
\juerg{M. Hinze and C. Kahle, A nonlinear model predictive concept for control of two-phase flows governed by the Cahn--Hilliard Navier--Stokes system. {\it System Modeling and Optimization}, 25th IFIP TC 7 Conference 2011, IFIP AICT 391: (2012) 348-357.}

\bibitem{Is} 
\pier{H.\ Israel,
Long time behavior of an {A}llen-{C}ahn type equation with a
singular potential and dynamic boundary conditions,
{\it J. Appl. Anal. Comput.} {\bf 2} (2012) 29-56.%
}

\bibitem{Lions}
J.-L. Lions,
``\'Equations diff\'erentielles op\'erationnelles et probl\`emes aux limites'',
Grundlehren, Band~111,
Springer-Verlag, Berlin, 1961.

\bibitem{LioMag}
J.-L. Lions and E. Magenes,
``Non-homogeneous boundary value problems and applications'',
Vol.~I,
Springer, Berlin, 1972.

\bibitem{MirZelik} 
A. Miranville and S. Zelik,
Robust exponential attractors for Cahn--Hilliard type equations with singular potentials,
{\it Math. Methods Appl. Sci.\/} 
{\bf 27} (2004) 545--582.

\bibitem{Podio}
P. Podio-Guidugli, 
Models of phase segregation and diffusion of atomic species on a lattice,
{\it Ric. Mat.} {\bf 55} (2006) 105-118.

\bibitem{PRZ} 
J. Pr\"uss, R. Racke and S. Zheng, 
Maximal regularity and asymptotic behavior of solutions for the Cahn--Hilliard equation with dynamic boundary conditions,  
{\it Ann. Mat. Pura Appl.~(4)\/}
{\bf 185} (2006) 627-648.

\bibitem{RZ} 
R. Racke and S. Zheng, 
The Cahn--Hilliard equation with dynamic boundary conditions, 
{\it Adv. Differential Equations\/} 
{\bf 8} (2003) 83-110.

\bibitem{RoSp}
\juerg{E. Rocca and J. Sprekels,
Optimal distributed control of a nonlocal convective Cahn--Hilliard equation 
by the velocity in 3D, WIAS Preprint No. 1942 (2014), pp. 1-24.} 

\bibitem{Simon}
J. Simon,
{Compact sets in the space $L^p(0,T; B)$},
{\it Ann. Mat. Pura Appl.~(4)\/} 
{\bf 146} (1987) 65-96.

\bibitem{WaNa}
Q.-F. Wang and S.-I. Nakagiri, Optimal control of distributed parameter 
system given by Cahn--Hilliard equation, {\it Nonlinear Funct. Anal. Appl.} 
{\bf 19} (2014) 19-33.

\bibitem{WZ} H. Wu and S. Zheng,
Convergence to equilibrium for the Cahn--Hilliard equation with dynamic boundary conditions, 
{\it J. Differential Equations\/}
{\bf 204} (2004) 511-531.

\bibitem{ZL1}
\juerg{X. P. Zhao and C. C. Liu, Optimal control of the convective Cahn--Hilliard equation, 
{\it Appl. Anal.\/} {\bf 92} (2013) 1028-1045.}

\bibitem{ZL2}
\juerg{X. P. Zhao and C. C. Liu, Optimal control for the convective Cahn--Hilliard equation in 2D
case, {\it Appl. Math. Optim.} \pier{{\bf 70} (2014) 61-82.}}

\End{thebibliography}

\End{document}

\bye